\documentclass[11pt]{article}
\usepackage{authblk}
\usepackage{amsmath,amstext,amsfonts,amssymb,amsthm} %
\usepackage{graphicx,color} %
\usepackage[latin1]{inputenc}   %
\usepackage[T1]{fontenc} %
\usepackage[english]{babel} %
\usepackage[caption=false,font=normalsize,labelfont=sf,textfont=sf]{subfig}
\usepackage{textcomp}
\usepackage{stfloats}
\usepackage{url}
\usepackage{verbatim}
\usepackage{cite}
\usepackage{bm}
\usepackage{bbm}
\usepackage{tikz}
\usetikzlibrary{positioning}
\usepackage{geometry}
\usepackage{xcolor}

\usepackage{mathtools}
\usepackage{bm}
\usepackage{dsfont}
\usepackage{caption}
\usepackage{fullpage} 
\usepackage{mathrsfs} 
\usepackage{xcolor}
\usepackage{marginnote} %

\usepackage[T1]{fontenc} 
\usepackage[english]{babel} 
\usepackage[plainpages=false,pdfpagelabels,colorlinks=true,linkcolor=red,urlcolor=red,citecolor=blue]{hyperref}

\let\OLDthebibliography\thebibliography
\renewcommand\thebibliography[1]{
  \OLDthebibliography{#1}
  \setlength{\parskip}{0pt}
  \setlength{\itemsep}{0pt plus 0.6ex}
}

\makeatletter

\usepackage{kantlipsum}

\allowdisplaybreaks

\numberwithin{equation}{section}
\newtheorem{theorem}{Theorem}[section]

\newtheorem{lemma}{Lemma}[section]

\newtheorem{proposition}{Proposition}[section]
\newtheorem{definition}{Definition}[section]

\newtheorem{remark}{Remark}[section]

\usepackage[plainpages=false,pdfpagelabels,colorlinks=true,linkcolor=red,urlcolor=red,citecolor=blue]{hyperref}

\usepackage[symbol]{footmisc}

\renewcommand{\thefootnote}{\fnsymbol{footnote}}
\def\ex{\mathbb{E}}

\def\exp{\text{exp}}
\def\X{\mathbf{x}}
\def\Y{\mathbf{Y}}
\def\Z{\mathbf{Z}}

\def\bU{\mathbf{U}}

\def\bA{\mathbf{A}}

\def\bQ{\mathbf{Q}}
\def\bS{\mathbf{S}}
\def\bI{\mathbf{I}}
\def\G{\mathbf{G}}
\def\bg{\mathbf{g}}
\def\bz{\mathbf{z}}
\def\bT{\mathbf{T}}

\def\Var{\mathrm{Var}}
\def\tr{\mathrm{Tr}}
\def\Tr{\mathrm{Tr}}

\def\E{\mathbb{E}}

\def\LLambda{\boldsymbol{\Lambda}}
\def\PPhi{\boldsymbol{\Phi}}
\def\pphi{\boldsymbol{\phi}}
\def\bbeta{\boldsymbol{\beta}}

\DeclareMathOperator*{\argmin}{arg\,min}

\title{\textbf{Inference in High-Dimensional Time Series Regression Models}}
\author{Samriddha Lahiry\textsuperscript{*}\textsuperscript{$\dagger$} ,
Daria Tieplova\textsuperscript{*}\textsuperscript{$\ddagger$},}

\date{}

\begin{document}

\maketitle

\begingroup
\renewcommand{\thefootnote}{\fnsymbol{footnote}}
\setcounter{footnote}{0}
\footnotetext{\textsuperscript{*}Authors are listed in alphabetical order of their surnames.}
\footnotetext{\textsuperscript{$\dagger$}Department of Statistics and Data Science, National University of Singapore, email: slahiry@nus.edu.sg}
\footnotetext{\textsuperscript{$\ddagger$}Department of Mathematics, Aarhus University \& B. Verkin Institute for Low Temperature Physics and Engineering  of the National Academy of Sciences of Ukraine, email: dtieplova@math.au.dk}
\endgroup

\begin{abstract}
  High-dimensional regression with temporally dependent observations arises in a variety of applications. While the sparse setting has been extensively studied, the dense-signal regime remains comparatively poorly understood. We study ridge regression under a row-dependent design in which different groups of covariates evolve according to distinct autoregressive processes, allowing the number of features and the sample size to grow at comparable rates without imposing sparsity on the underlying signal. Using tools from random matrix theory, we derive a precise characterization of the asymptotic estimation risk. We then address the bias induced by ridge regularization by constructing a debiased estimator and establishing its joint asymptotic normality for any fixed collection of coordinates. We characterize the corresponding limiting variances and estimate them through joint estimation of the blockwise signal energies and the observation-noise variance, thereby enabling valid statistical inference. Furthermore, we establish a correspondence between our deterministic-equivalent equations and the Bayesian fixed-point equations under a Gaussian prior, providing a direct proof of the Gaussian specialization of a previously conjectured MMSE formula. Our results thus provide a frequentist counterpart to recent Bayesian and information-theoretic analyses of the same model and extend the proportional asymptotic theory of regularized regression beyond row-independent designs.
\end{abstract}
\section{Introduction}

High-dimensional data exhibiting both cross-sectional and temporal dependence are ubiquitous, with applications in macroeconomic forecasting \cite{li2014forecasting}, genomics \cite{bar2012studying}, systemic risk measurement in financial markets \cite{basu2024high}, and fMRI studies \cite{sporns2021dynamic}. Beyond the challenges posed by high dimensionality, dependence among observations introduces additional complications and calls for methods tailored to the underlying dependence structure. Over the past decade, techniques from high-dimensional statistics have been extended to such dependent settings, most notably in the seminal work of \cite{basu2015regularized}. Existing results, however, have largely relied on low-complexity structural assumptions, typically requiring the regression signal to be highly sparse. Although these assumptions facilitate mathematical analysis, they can be restrictive in practice and limit the applicability of existing methods when the underlying signal is not sufficiently sparse.

In parallel, a new paradigm has emerged in high-dimensional asymptotics, motivated primarily by ideas from statistical physics and information theory. This framework relaxes stringent sparsity assumptions and enables the analysis of models in the proportional asymptotic regime, where the number of features $p$ and the number of observations $N$ grow at comparable rates. In the frequentist setting, a broad toolkit has been developed for signal estimation in this regime. Prominent approaches include approximate message passing
\cite{DMM09,bayati_mont,Sur_2,lnk_flra_survey,AMP_survey}; leave-one-out arguments
\cite{elkaroui_pnas_2,elkaroui_ptrf,chen2021spectral}; and the convex Gaussian min-max theorem
\cite{stojnic_1,stojnic_2,CT_5}. A central advantage of this perspective is that it dispenses with stringent sparsity requirements and naturally accommodates structured signal classes that may be either dense or sparse. Bayesian counterparts to these techniques have also been developed. In particular, linear regression with an i.i.d.\ Gaussian design has been studied extensively \cite{barbier2016mutual,barbier2018mutual,reeves_2016,barbier2019optimal}, although extensions to more structured designs remain limited.

Several structured designs arising from dependent random-matrix ensembles have recently been analyzed in the proportional asymptotic regime \cite{dudeja2022universality,dudeja2023universality,dudeja2024spectral,lahiry2024universality,saeed_mont}. These ensembles, however, do not capture the temporal dependence across rows that naturally arises in time-series models. Motivated by \cite{basu2015regularized}, the recent work \cite{mutinf_time} considers a stochastic regression model in which the rows of the design follow an AR(1) process. Whereas \cite{basu2015regularized} studies this model under sparsity assumptions, \cite{mutinf_time} analyzes it in the proportional asymptotic regime without imposing sparsity. The latter work adopts a Bayesian formulation in which a prior is placed on the signal and seeks to characterize the minimum mean-squared error (MMSE) of signal estimation, a quantity closely connected to the mutual information between the signal and the observations. More specifically, its objectives are to derive an asymptotic formula for the mutual information per observation as $N\to\infty$, characterize the limiting MMSE, and construct a computationally efficient estimator that asymptotically attains the MMSE.

In this paper, we study the model introduced in \cite{mutinf_time} from a frequentist perspective, thereby providing a complementary counterpart to its Bayesian analysis. We consider ridge-penalized estimation under row-dependent designs and characterize the estimator's asymptotic behavior in the proportional regime. Unlike in classical sparse regimes, estimation in the dense proportional regime is generally inconsistent: the estimation risk converges to a nonzero constant characterized by the solution of a fixed-point equation. Although this phenomenon is well understood for row-independent designs, we establish an analogous characterization for ridge regression under the row-dependent model considered here. We then address the bias inherent in ridge-penalized estimation. Using tools from random matrix theory, we construct a debiased ridge estimator and establish its asymptotic normality. Further we compute the explicit formula for the limiting variance and construct its estimator thereby enabling valid coordinate-wise inference for the underlying signal.

The scope of these results extends beyond the specific AR(1) design studied here. Temporal dependence across observations is a recurring feature of financial econometric models, particularly high-dimensional predictive regressions \cite{mei2024lasso,gao2024econometric}. It also arises naturally in vector autoregressive (VAR) models, where each component of a multivariate time series is regressed on lagged values of the full system and the resulting lead--lag relationships are central objects of interest; see, for example, \cite{basu2015regularized}. Although these models involve additional features not captured by our setting, such as endogenous lagged covariates and more intricate cross-sectional dynamics, they share the fundamental difficulty that dependence across rows alters the high-dimensional behavior of the sample covariance matrix and, consequently, that of regularized estimators.
Closely related challenges arise in spatial and spatio-temporal models. Spatial observations typically exhibit dependence determined by geographic or network proximity \cite{cressie2015statistics}, and recent work has begun to address regression and variable selection in high-dimensional spatial settings \cite{hou2024estimation,pfarrhofer2019flexible}. Our results do not directly encompass these models, whose dependence structures can be substantially more complex than the temporal mechanism considered here. Nevertheless, even for the comparatively simple row-dependent design studied in this paper, obtaining a precise characterization of ridge regression and constructing a valid debiasing procedure require a delicate analysis of the associated dependent random matrices. To the best of our knowledge, such questions remain largely unexplored (cf. \cite{Li2026} for a specific case). We therefore view the present work as an initial step toward a unified high-dimensional theory for dependent data-one that could eventually accommodate temporal, spatial, and spatio-temporal structures and address both frequentist and information-theoretic questions. Developing this broader framework is an important direction for future research.

\subsection{Major contributions}

We study estimation and inference in high-dimensional linear regression
with temporally dependent covariates in the proportional asymptotic
regime, without imposing sparsity assumptions on the underlying signal.

\paragraph{1. Characterization of the asymptotic estimation risk.}

We first derive the exact asymptotic estimation risk of the
ridge-penalized estimator. As is typical in the proportional asymptotic
regime, the limiting risk is characterized through the solution of a
system of fixed-point equations. The risk formula accounts for
heterogeneous temporal covariance structures and different signal
energies across blocks, and is expressed in terms of deterministic
equivalents and their derivatives with respect to the regularization
parameter. For the autoregressive model, the Toeplitz structure of the
covariance matrices yields explicit integral representations in terms
of the corresponding spectral densities.

\paragraph{2. Debiasing and joint asymptotic normality.}

Our second and principal contribution is to establish the joint
asymptotic normality of a debiased version of the ridge estimator.
Using deterministic equivalents, we identify block-specific debiasing
factors that remove the leading regularization bias.  The  joint CLT of the debiased ridge estimator has a diagonal limiting
covariance matrix, with explicit block-specific variances that depend
on the observation-noise variance and the signal energies across all
blocks. The argument
extends beyond the autoregressive setting to independent Gaussian
columns with more general block covariance matrices, under the
conditions of our theorem.

\paragraph{3. Consistent variance estimation and feasible inference.}

A particular difficulty in the proportional asymptotic regime is that
the limiting variances depend on unknown characteristics of the
high-dimensional signal, even when the design covariance matrices are
known. Motivated by the method-of-moments approach of
\cite{dicker2014variance}, we construct observable quadratic moments
that jointly identify the blockwise signal energies and the
observation-noise variance.
We derive exact finite-dimensional moment
identities that account for the dependence between the response and
the block sample covariance matrices, and establish uniform
invertibility of the resulting system.
The estimators attain the parametric rate relative to their
finite-dimensional targets and are independent of the choice of the
ridge penalty. Substituting them into the asymptotic variance formula
yields consistent variance estimators and, through the joint CLT,
asymptotically valid confidence intervals, joint confidence regions
and hypothesis tests for finitely many coefficients.

\paragraph{4. Connection with Bayesian estimation.}

The Bayesian analysis of the same time-series model in
\cite{mutinf_time} establishes limiting formulas for the mutual
information and measurement MMSE, and conjectures a formula for the
blockwise parameter MMSE. Under a centered Gaussian prior, ridge with
the matched penalty coincides with the posterior mean. We show that
our deterministic-equivalent equations then coincide, after a change
of variables, with the Bayesian fixed-point equations. This proves
the Gaussian specialization of the blockwise MMSE formula and recovers
the measurement MMSE and mutual-information formulas independently of
the compact-support assumption in that work.

\subsection{Related literature}

Over the past decade, substantial progress has been made in developing
valid inferential procedures based on debiased penalized regression
estimators \cite{van2014asymptotically,zhang_zhang,javanmard2014confidence}.
Much of this literature relies on sparsity or additional structural
assumptions on the design covariance. Our work is related to the
leave-one-out analysis of \cite{elkaroui_ptrf}: both approaches exploit
resolvent representations, but their setting permits dependence across
covariates while assuming independent observations, whereas ours
allows dependence across observations through independent Gaussian
columns with heterogeneous temporal covariances.

Our debiasing principle is closely aligned with that of \cite{deb_med},
whose analysis relies on the bilinear-form fluctuation results of
\cite{yin2023limiting}; see also \cite{bai_eigen} for related
eigenvector asymptotics. To illustrate the distinction between the
two settings, let $G\in\mathbb R^{N\times p}$ have independent standard
Gaussian entries. For independent observations with covariate covariance
$\Sigma$ and for observations with common temporal covariance $\Lambda$,
respectively, the relevant resolvents are
\[
  \left(p^{-1}\Sigma^{1/2}G^\top G\Sigma^{1/2}
        +\kappa I_p\right)^{-1},
  \qquad
  \left(p^{-1}G^\top\Lambda G+\kappa I_p\right)^{-1}.
\]
Although the spectral relationships and deterministic equivalents for
the latter are well understood through the spectral properties of
the companion matrix $p^{-1}\Lambda^{1/2}GG^\top\Lambda^{1/2}$
\cite{Couillet2011,benaych2016spectral}, fluctuations of bilinear forms
depend additionally on eigenvectors and cannot, in general, be
transferred through spectral identities alone. A key technical
contribution of this work is the characterization of the fluctuations
of bilinear forms of
$\left(p^{-1}G^\top\Lambda G+\kappa I_p\right)^{-1}$,
which enables successful debiasing.

A related approach is the spectrum-aware debiasing framework of
\cite{Li2026}, which accommodates non-sparse signals and
right-rotationally invariant designs. The authors exploit this invariance,
together with the $R$-transform, to establish a central limit theorem for
debiased estimators. The single-block case of our model is a special case
of their framework. However, their methods do not directly apply to the general
blockwise right-rotationally invariant setting considered here. 

Finally we note that our block model also arises in multiuser communication systems,
where \cite{cakmak2026} studies an orthogonal approximate message passing
algorithm and relates its asymptotic performance to the predicted
Bayes-optimal risk. A related analysis in the Bayesian framework has been considered in \cite{mutinf_time}. These algorithmic and Bayesian investigations are
complementary to our frequentist analysis.
\subsection{Outline of the paper}

In Section~\ref{sec:prelim}, we introduce the ridge estimator and briefly
describe the heuristics underlying the construction of the debiased
estimators. Section~\ref{sec:main} presents the main results. In
particular, Theorem~\ref{thm:asymp_risk} characterizes the asymptotic
risk of the ridge estimator, while Theorem~\ref{thm:joint_clt}
establishes the asymptotic normality of the debiased ridge estimator.
The corresponding asymptotic variances are characterized through the
solution of a system of equations. In
Section~\ref{sec:consistent_variance}, we construct consistent
estimators of these limiting variances.
Theorem~\ref{thm:var_est_studentized} establishes the asymptotic
normality of the resulting studentized estimator, thereby enabling
valid inference on the regression coefficients. In
Section~\ref{sec:bayes_conn}, we discuss connections with the Bayesian
counterpart of the model, while
Section~\ref{sec:ridge_simulations} presents numerical experiments
supporting the asymptotic normality of the proposed debiased estimator.
Finally, Section~\ref{sec:discussion} concludes with a discussion of
several future direction.

\subsection{Notations}
 For all vectors and matrices we will be using bold notations $\Y_N,\boldsymbol{\Phi},$ etc. and their elements will be $Y_{N,i}, \Phi_{ij},$ etc. $\mathbf{I}_n$ represents the identity matrix of size $n\times n$. We use the standard notations $O(.)$ and $o(.)$ to denote ``asymptotically bounded'' and ``asymptotically vanishing'' respectively, i.e., we say $a_N=O(\gamma_N)$ if $\limsup|a_N/\gamma_N|\leq C$ and $a_N=o(\gamma_N)$ if $\limsup|a_N/\gamma_N|=0$. 
 Expectation with respect to a group of random variables (conditioned on the rest of the random variables) will be denoted as $\ex_X[.]$  where $X$ would denote the relevant collection of random variables whereas  $\ex[.]$ and $\langle .\rangle$ will be reserved for unconditional expectation. Norm $\|.\|$ is the $\ell_2$ norm of a vector, unless otherwise mentioned. For matrices the $\|\bA\|$ will always denote the largest singular value of $\bA$ (operator norm) and hence the same notation $\|.\|$ will be reserved for the $\ell_2$ norm of a vector and the operator norm of a matrix, the difference being understood from the context. Throughout the notation $C$ will denote a generic constant unless otherwise specified, that may change from place to place. The Kronecker delta function will be denoted as  
$\delta_{\mu,\nu}=1_{\mu=\nu}$.

\section{Preliminaries}\label{sec:prelim}

Let $(Y_\mu,\X_\mu)_{\mu=1}^N$ be $N$ observations from a time
series, with the dependence structure outlined below. Let $\X_\mu$
be a $p$-dimensional, centered, stationary Gaussian process given by
\begin{equation}\label{eq:def_x_t}
    \X_{\mu+1}=\bA_p\X_\mu+\boldsymbol{\xi}_{\mu},
    \quad \mu\in\mathbb{Z},
\end{equation}
where $\bA_p$ is a deterministic $p\times p$ matrix and
$(\boldsymbol{\xi}_\mu)_{\mu\in\mathbb Z}$ are i.i.d.\
$\mathcal N(0,\sigma_1^2\bI_p)$ innovations, with $\sigma_1>0$.
For every $\mu$, $\boldsymbol{\xi}_\mu$ is independent of
$\sigma(\X_\nu:\nu\leq\mu)$.
The observed time series $Y_\mu$ has the form
\begin{equation}\label{eq_SRM}
    Y_\mu=\dfrac{1}{\sqrt{p}}\X_\mu^\intercal\bbeta_0+Z_\mu,
    \quad \mu=1,\ldots,N,
\end{equation}
where $\bbeta_0$ denotes the true signal vector and the noise variables
$Z_\mu$ are i.i.d.\ Gaussian with mean zero and variance $\sigma_2^2$,
independent of the design process. In our frequentist analysis, we write
$\bbeta=\bbeta_0$ and treat $\bbeta$ as deterministic.
By rescaling both $\X_\mu$ and $Y_\mu$ by $\sigma_1^{-1}$,
we obtain an equivalent model with unit innovation variance and
observation-noise variance $\sigma_2^2/\sigma_1^2$, without changing
the signal vector. Thus, without loss of generality, we work with
unit innovation variance and denote the resulting observation-noise
variance by $\sigma^2$. 

We are interested in the high-dimensional regime where the number of observations $N$ and the dimension $p$ both tend to infinity at the same rate. That is, we assume that $p=p(N)$ and for $c_p:=N/p$ we have
$$\lim_{N\rightarrow \infty}c_p =\lim_{N\to \infty}\frac Np=c \in (0,\infty).$$ 
 Let $\boldsymbol{\Phi}_N=(\X_1^\intercal,\ldots, \X_N^\intercal)^\intercal$ be a $N\times p$ matrix whose rows are the vectors $(\X_\mu)_{\mu=1}^N$. Further, let $\Y_N=(Y_1,\ldots, Y_N)^\intercal$ and $\Z_N=(Z_1,\ldots,Z_N)^\intercal$ be two vectors of length $N$ containing the individual observations and noises, respectively. Then our model can be written in a vector form:
\begin{align}\label{eq:mod_vec}
    \Y_N=\dfrac{1}{\sqrt{p}}\boldsymbol{\Phi}_N\bbeta_0+\Z_N.
\end{align}

The form of $\bA_p$ plays a major role in the type of resulting linear regression. An important example, addressed in, e.g., \cite{Takeda_2006,tulino2013support,barbier2018mutual,rangan2019vector} (with the final work proposing a near-optimal algorithm for the estimation of signal $\boldsymbol{\beta}_0$), is when $\bA_p$ is a symmetric random matrix with uniformly distributed eigenvectors. Then, one has the eigendecomposition $\bA_p=\boldsymbol{U}\boldsymbol{D}\boldsymbol{U}^{-1}$ with $\boldsymbol{U}$ a Haar distributed random matrix and we may rewrite \eqref{eq:def_x_t} as $\boldsymbol{U}^{-1}\X_{\mu+1}=\boldsymbol{D}\boldsymbol{U}^{-1}\X_\mu+\boldsymbol{\xi}_{\mu}$, noting that $\boldsymbol{U}^{-1}\boldsymbol{\xi}_{\mu}=\boldsymbol{\xi}_{\mu}$ in law due to rotational invariance of the Gaussian noise. Thus the time series $\bz_\mu=\boldsymbol{U}^{-1}\X_\mu$ satisfies $\bz_{\mu+1}=\boldsymbol{D}\bz_\mu+\boldsymbol{\xi}_{\mu}$, which upon rewriting $\boldsymbol{\Phi}_N
=\tilde{\boldsymbol{\Phi}}_N\boldsymbol{U}^{\intercal}$ with $\tilde{\boldsymbol{\Phi}}_N=(\bz_1^\intercal,\ldots,\bz_N^\intercal)^\intercal$ brings us to the linear regression model with a right-rotationally invariant design matrix. 

Our present interest lies in the extreme opposite setting where the previously Haar matrix $\boldsymbol{U}$ is now proportional to the identity so that $\bA_p$ becomes diagonal, with a rather general choice of eigenvalues. Such a case \emph{cannot} be treated with the aforementioned line of work and requires new techniques, both from the analytical and algorithmic perspective. We thus denote
$\bA_p={\rm diag}(\lambda_1,\lambda_2,\ldots, \lambda_p)$, where $\lambda_i$'s assumed to be known parameters, s.t. $|\lambda_i|<1$ for all $i$. Under this assumption, equation (\ref{eq:def_x_t}) breaks down into $p$ independent equations $$x_{\mu+1,i}=\lambda_i x_{\mu,i}+\xi_{\mu,i},\quad i=1,\ldots,p\,,$$ 
which allows us to easily calculate the variance of $x_{\mu,i}$ and the correlations between $x_{\mu,i}$ and $x_{\mu+s,i}$:
\begin{align}
    \ex(x_{\mu+1,i})^2=\lambda_i^2\ex(x_{\mu,i})^2+1.
\end{align}
Since $\X_\mu$ is a stationary process, its variance does not depend on $\mu$ and the last equality gives us:
\begin{align}\label{eq:var_x_t^i}
    \ex(x_{\mu,i})^2=(1-\lambda_i^2)^{-1}.
\end{align}
From this we immediately obtain 
\begin{align}\label{eq:cor_x_t^i}
    \ex(x_{\mu,i}x_{\mu+s,i})=\ex(x_{\mu,i}(\lambda_i^sx_{\mu,i}+\sum_{j=0}^{s-1}\lambda_i^j\xi_{\mu-j+s-1,i}))=(1-\lambda^2_i)^{-1}\lambda_i^s.
\end{align}

Independence of coordinates also implies that the columns of $\boldsymbol{\Phi}_N$ are independent Gaussian vectors and, due to (\ref{eq:var_x_t^i}) and (\ref{eq:cor_x_t^i}), 
we can easily calculate the covariance matrix $\LLambda_{i,N}$ of the $i$-th column $(x_{1,i},x_{2,i},\ldots,x_{N,i})^\intercal$:
\begin{align}\label{eq:def_Lambda}
    \LLambda_{i,N}:=\dfrac{1}{1-\lambda_i^2}\begin{pmatrix}
    1&\lambda_i &\ldots&\lambda_i^{N-1}\\
    \lambda_i&1&\ldots&\lambda_i^{N-2}\\
    \vdots&\vdots&\ddots&\vdots\\
    \lambda_i^{N-1}&\lambda_i^{N-2}&\ldots&1
\end{pmatrix}.
\end{align}
In what follows, we will often drop indices $N,p$ when it is not crucial. 
With  this, we can present the design matrix as $\boldsymbol{\Phi}_N=(\LLambda_{1,N}^{1/2}\bg_{1},\ldots,\LLambda_{p,N}^{1/2}\bg_{p})$, where the $\bg_i$ are i.i.d. standard Gaussian vectors. In this form it is easy to see the connection to the previously well-studied models.

\textit{Case 1: $\lambda_i=0$ for $i=1,\ldots,p.$} 
In this case all $\LLambda_{i,N}$ are the identity and our model collapses to that of standard linear regression, where the design matrix has i.i.d. standard Gaussian elements. This model was studied in the CDMA \cite{tanaka2002statistical}
 and compressed sensing contexts \cite{DMM09,krzakala2012statistical,barbier2016mutual,barbier2020mutual,reeves_2016}. 

\textit{Case 2: $\lambda_i=\lambda\neq0$ for $i=1,\ldots,p.$} 
Again, all columns of $\boldsymbol{\Phi}_N$ will share the same covariance matrix, denoted by  $\LLambda$. This allows us to rewrite $\boldsymbol{\Phi}_N=\LLambda^{1/2}\G_N$, where $\G_N$ is an $N\times p$ matrix with i.i.d. standard Gaussian elements. Since $\G_N$ is rotationally invariant, this results in the design matrix $\boldsymbol{\Phi}_N$  belonging to the class of right-rotationally invariant matrices. Single-letter formulas for the mutual information of such models was obtained in \cite{Takeda_2006,tulino2013support,barbier2018mutual}, while CLT for a debiased estimator has been obtained in \cite{Li2026}.

In this work, we consider the more general case in which the diagonal entries
of $\bA_p$ take $k$ distinct values, with $k$ fixed as $p\to\infty$.
This is a natural first step toward studying the model with an arbitrary diagonal
matrix $\bA_p$. 

\textit{Case 3:}\label{Assum:A} 
We assume that diagonal matrix $\bA_p$ has $k$ different eigenvalues $\lambda_1,\ldots,\lambda_k$. Further, we assume that $k$ and the distinct known parameters
$\lambda_1,\ldots,\lambda_k\in(-1,1)$ are fixed as
$N,p\to\infty$. There is a partition
$\{1,\ldots,p\}=\bigcup_{a=1}^k I_a$ such that
$A_{jj}=\lambda_a$ for $j\in I_a$. Importantly, we require that the cardinality of each such set grows proportionally with $p$, i.e.,  we have $|I_i|/p\rightarrow \ell_i>0$ for each $i=1,\ldots,k$.

Without loss of generality, we further assume that  $I_1$ contains the first $|I_1|$ integers, $I_2$ contains the next $|I_2|$ consecutive integers and so on, since it affects only the permutation of the columns of $\boldsymbol{\Phi}_N$. Then,
 $\boldsymbol{\Phi}_N$ can be represented as a block matrix with $k$ blocks  of dimensions $N\times |I_i|$, i.e., $(\LLambda_{1,N}^{1/2}\mathbf{G}_N^{(1)},\ldots,\LLambda_{k,N}^{1/2}\mathbf{G}_N^{(k)})$, where for every $i=1,\ldots,k$,  $\mathbf{G}_N^{(i)}$ is an $N\times |I_i|$ independent matrix with i.i.d. standard Gaussian elements. Thus (\ref{eq:mod_vec}) can be rewritten as
\begin{align}
    \Y_N=\frac{1}{\sqrt{p}}(\LLambda_{1,N}^{1/2}\mathbf{G}_N^{(1)},\ldots,\LLambda_{k,N}^{1/2}\mathbf{G}_N^{(k)})\bbeta_0+\Z_N \label{diagonal_lambda_model}.
\end{align}
The block structure of the design matrix naturally propagates to the signal. In what follows we define   $\bbeta_0^{(i)} := (\beta_{0,j})_{j\in I_i} \in \mathbb{R}^{|I_i|}$ as the $i$-th block of $\bbeta_0$, and similarly for $\bbeta^{(i)}$.
We emphasize that in this setting the design matrix is not right-rotationally invariant, contrary to the Cases 1 and 2 described above. Instead, it exhibits a richer "block-wise" right-rotational invariance. 

In this paper, we consider a penalized estimation of the signal $\bbeta_0$. In particular we define our estimator as follows:
$$\hat{\bbeta}(\kappa)=\argmin_{\bbeta}\frac{1}{2}\|\Y_N-\dfrac{1}{\sqrt{p}}\boldsymbol{\Phi}_N\bbeta\|^2+\frac{\kappa}{2} \|\bbeta\|^2$$
which can be further simplified to
\begin{equation}{\label{eqn:ridge_closed}}\hat{\bbeta}(\kappa)=p^{-1/2}(\Phi_N^\intercal\Phi_N/p+\kappa I)^{-1}\Phi_N^\intercal\Y_N.\end{equation}
In the next section, we show how the theory of deterministic equivalent can be used to compute the asymptotic risk of the ridge estimator and show asymptotic normality of estimator formed by debiasing the ridge estimator. We describe the heuristics of the analysis below.

Using \eqref{eq:mod_vec}, we can write the ridge estimator as
\begin{align*}
\hat{\bbeta}(\kappa)
&=p^{-1/2}(\Phi_N^\intercal\Phi_N/p+\kappa I)^{-1}\Phi_N^\intercal
  \left(p^{-1/2}\Phi_N\bbeta+\Z_N\right)\\
&=\bbeta-\kappa(\Phi_N^\intercal\Phi_N/p+\kappa I)^{-1}\bbeta
  +p^{-1/2}(\Phi_N^\intercal\Phi_N/p+\kappa I)^{-1}\Phi_N^\intercal\Z_N.
\end{align*}
From this expression, we can compute the asymptotic risk
$p^{-1}\|\hat{\bbeta}(\kappa)-\bbeta\|^2$ of the ridge estimator.

It is known that $(\Phi_N^\intercal\Phi_N/p+\kappa I)^{-1}$ has a block diagonal deterministic equivalent (see the next section). For $j\in I_a$, let $\breve\beta_j=\mathcal C_a\hat\beta_j$,
where the factor $\mathcal C_a$ is given in \eqref{eq:debiased_ridge}.
The above identity then gives
\begin{align*}
\breve\beta_j-\beta_j
={}&\left[(\mathcal C_a-1)\beta_j
  -\kappa\mathcal C_a\mathbf e_j^\intercal
   (\Phi_N^\intercal\Phi_N/p+\kappa I)^{-1}\bbeta\right]\\
&+\mathcal C_a p^{-1/2}\mathbf e_j^\intercal
  (\Phi_N^\intercal\Phi_N/p+\kappa I)^{-1}\Phi_N^\intercal\Z_N.
\end{align*}
We will show CLT for the bracketed term by directly analyzing the fluctuation of the bilinear form of the resolvent. The noise term in the second line is conditionally Gaussian and also contributes to the limiting variance.
 Under the normalization in
\eqref{eq:mod_vec}, Theorem~\ref{thm:joint_clt} gives a joint Gaussian
limit for the debiased coordinate errors without any additional $\sqrt N$
normalization.

\section{Main results}\label{sec:main}

In this section, we state the main theorems of our paper.
We first introduce the resolvent and its companion
\begin{align}
    \bQ_p(\kappa)
    &:=
    \left(\frac1p\PPhi^\top\PPhi+\kappa\bI_p\right)^{-1},
    \label{eq:gram_resolvent}\\
    \tilde\bQ_p(\kappa)
    &:=
    \left(\frac1p\PPhi\PPhi^\top+\kappa\bI_N\right)^{-1},
    \label{eq:companion_resolvent}
\end{align}
where $\kappa>0$.

Our results below rely on the theory of deterministic-equivalents, which we recall here briefly. By the notation $A_p\longleftrightarrow B_p$ we will mean the asymptotic equivalence of the two matrices $A_p$ and $B_p$ in the following sense:
\[
    \frac1p\Tr\!\left[D_p(A_p-B_p)\right]
    \xrightarrow{\mathbb P}0,
    \qquad
    \mathbf u_p^\top(A_p-B_p)\mathbf w_p
    \xrightarrow{\mathbb P}0,
\]
for all deterministic matrices $D_p$ and vectors
$\mathbf u_p,\mathbf w_p$ of compatible dimensions
with uniformly bounded norms. If $B_p$ is deterministic we will call $B_p$ deterministic equivalent of $A_p$ (see Definition \ref{def_det_eq} for further details).

The following deterministic equivalent identities from
\cite{benaych2016spectral}, adapted to our normalization, will be crucial for establishing the risk asymptotics and the debiased CLT:
\begin{align}
    \bQ_p(\kappa)
    &\longleftrightarrow
    \bT_p(\kappa)
    :=
    \operatorname{diag}\left\{
        m_a(\kappa)\bI_{|I_a|}
    \right\}_{a=1}^k,
    \label{eq:gram_det_eq}\\
    \tilde\bQ_p(\kappa)
    &\longleftrightarrow
    \tilde\bT_p(\kappa)
    :=
    \frac1\kappa
    \left(
        \bI_N+\sum_{b=1}^k
        \frac{|I_b|}{p}m_b(\kappa)\Lambda_b
    \right)^{-1}.
    \label{eq:companion_det_eq}
\end{align}
Here, $\{m_a(\kappa),\tilde m_a(\kappa)\}_{a=1}^k$
are the unique positive solutions of
\begin{align}
    m_a(\kappa)
    &=\frac1{\kappa+\tilde m_a(\kappa)}, \label{fp_1}\\
    \tilde m_a(\kappa)
    &=\frac1p\Tr\!\left[
        \Lambda_a
        \left(
            \bI_N+\sum_{b=1}^k
            \frac{|I_b|}{p}m_b(\kappa)\Lambda_b
        \right)^{-1}
    \right]. \label{fp_2}
\end{align}
 We write $m_a'(\kappa)$ for the derivative
with respect to $\kappa$.

\subsection{Characterizing the asymptotic risk}

\begin{theorem}\label{thm:asymp_risk}
Consider the block model \eqref{diagonal_lambda_model}. Suppose that $\bbeta$ is deterministic and that $\|\bbeta\|^2/p=O(1)$. For fixed $\kappa>0$, the risk of the ridge estimator satisfies
\begin{align}
    \frac{1}{p}\|\hat{\bbeta}-\bbeta\|^2
    =    -\kappa^2\sum_{a=1}^k    \frac{\|\bbeta^{(a)}\|^2}{p}    m_a'(\kappa)
    +    \sigma^2\sum_{a=1}^k    \frac{|I_a|}{p}    \left(m_a(\kappa)+\kappa m_a'(\kappa)\right)
    +o_p(1).
\end{align}
\end{theorem}

\begin{proof}[Proof of Theorem \ref{thm:asymp_risk}]
To characterize the asymptotic risk note that
\begin{align}
\hat{\bbeta}
=p^{-1/2}\bQ_p(\kappa)\Phi_N^\intercal
\left(\frac1{\sqrt p}\Phi_N\bbeta+\Z_N\right)
=\bbeta-\kappa\bQ_p(\kappa)\bbeta
+p^{-1/2}\bQ_p(\kappa)\Phi_N^\intercal\Z_N.
\label{ridge_decomp}
\end{align}
so that
\begin{align}
\|\hat{\bbeta}-\bbeta\|^2/p
=&p^{-1}\left\|-\kappa\bQ_p(\kappa)\bbeta+p^{-1/2}\bQ_p(\kappa)\Phi_N^\intercal\Z_N\right\|^2\nonumber\\
=&p^{-1}\kappa^2\bbeta^\intercal\bQ_p(\kappa)^2\bbeta+p^{-2}\Z_N^\intercal\Phi_N\bQ_p(\kappa)^2\Phi_N^\intercal\Z_N
-2\kappa p^{-3/2}\bbeta^\intercal\bQ_p(\kappa)^2
\Phi_N^\intercal\Z_N \nonumber\\
=&p^{-1}\kappa^2\bbeta^\intercal\bQ_p(\kappa)^2\bbeta
+p^{-2}\Z_N^\intercal\Phi_N\bQ_p(\kappa)^2
\Phi_N^\intercal\Z_N+o_p(1). \label{sq_risk_expansion}
\end{align}
where the last equality follows by argument given below.
Note that the conditional expectation of
$p^{-3/2}\bbeta^\intercal\bQ_p(\kappa)^2
\Phi_N^\intercal\Z_N$ given $\Phi_N$ is $0$,
and the variance is given by:
\[
\frac{\sigma^2}{p^2}
\bbeta^\intercal\bQ_p(\kappa)^2
\left(\frac{\Phi_N^\intercal\Phi_N}{p}\right)
\bQ_p(\kappa)^2\bbeta.
\]
Since
$\|\bQ_p(\kappa)^2(\Phi_N^\intercal\Phi_N/p)
\bQ_p(\kappa)^2\|=O_p(1)$,
and $\|\bbeta\|^2/p=O(1)$, the conditional variance is
$o_p(1)$ and hence the equality follows.

Now  we can consider the noise quadratic form. Applying
Lemma~\ref{quad_form_lemma} conditionally on $\Phi_N$ gives
\begin{equation}
\frac1{p^2}\Z_N^\intercal\Phi_N\bQ_p(\kappa)^2
\Phi_N^\intercal\Z_N
=\frac{\sigma^2}{p^2}\Tr\left[
\Phi_N\bQ_p(\kappa)^2\Phi_N^\intercal
\right]+o_p(1)
=\frac{\sigma^2}{p}\Tr\left[\bQ_p(\kappa)\right]
-\frac{\kappa\sigma^2}{p}\Tr\left[\bQ_p(\kappa)^2\right]
+o_p(1).\label{noise_quad}
\end{equation}
Lemma~\ref{quad_form_lemma} applies conditionally on $\Phi_N$, since
$\Z_N$ is independent of $\Phi_N$, and 
\[
\left\|
\frac1p\Phi_N\bQ_p(\kappa)^2\Phi_N^\intercal
\right\|\leq\frac1{4\kappa}.
\]

Next, we need the following lemma which deals with the square of resolvents.

\begin{lemma}\label{det_eq_quad}
For deterministic unit vectors $\alpha$ and $\gamma$, we have
\begin{equation}
\alpha^\intercal\bQ_p(\kappa)^2\gamma
=-\sum_{a=1}^k
\alpha^{(a)\intercal}\gamma^{(a)}m_a'(\kappa)+o_p(1). \label{bilin_quad}
\end{equation}
Moreover,
\begin{equation}
\frac1p\Tr\left[\bQ_p(\kappa)^2\right]
=-\sum_{a=1}^k\frac{|I_a|}{p}m_a'(\kappa)+o_p(1).\label{trace_quad}
\end{equation}
\end{lemma}
For convenience, we postpone the proof to the Appendix.
Combining \eqref{trace_quad} with the deterministic equivalents for $\bQ_p$ we have 
\begin{equation}\frac{\sigma^2}{p}\Tr\left[\bQ_p(\kappa)\right]
-\frac{\kappa\sigma^2}{p}\Tr\left[\bQ_p(\kappa)^2\right]
=\sigma^2\sum_{a=1}^k\frac{|I_a|}{p}m_a(\kappa)
+\kappa\sigma^2\sum_{a=1}^k\frac{|I_a|}{p}m_a'(\kappa)
+o_p(1).\label{trace_lin_quad}\end{equation}

Finally by applying \eqref{bilin_quad} to the normalized signal
$\bbeta/\|\bbeta\|$, we obtain
\begin{equation}
\frac1p\bbeta^\intercal\bQ_p(\kappa)^2\bbeta
=-\sum_{a=1}^k\frac{\|\bbeta^{(a)}\|^2}{p}
m_a'(\kappa)+o_p(1).
\label{qform_quad}
\end{equation}

If $\bbeta=0$, this identity is immediate. Combining \eqref{sq_risk_expansion}, \eqref{noise_quad}, \eqref{trace_lin_quad}, and \eqref{qform_quad} we obtain
\begin{align*}
\frac1p\|\hat{\bbeta}-\bbeta\|^2
=-\kappa^2\sum_{a=1}^k
\frac{\|\bbeta^{(a)}\|^2}{p}m_a'(\kappa)
+\sigma^2\sum_{a=1}^k\frac{|I_a|}{p}m_a(\kappa)+\kappa\sigma^2\sum_{a=1}^k\frac{|I_a|}{p}m_a'(\kappa)
+o_p(1),
\end{align*}
which proves the theorem.

\end{proof}

\subsection{Joint CLT}
For a coordinate $j\in I_a$, define the debiasing coefficient 
\begin{equation}\label{eq:debiased_ridge}
    \mathcal{C}_a    :=    (1-\kappa m_a(\kappa))^{-1},
\end{equation}
and the additional $k\times k$ matrix
\[
\Omega_{ab}(\kappa)
:=\ell_b\kappa^2m_b(\kappa)^2\frac1p\Tr\!\left(\Lambda_a\tilde\bT_p(\kappa)\Lambda_b\tilde\bT_p(\kappa)\right).
\]
\begin{theorem}\label{thm:joint_clt}
Consider the block model \eqref{diagonal_lambda_model}.
Suppose that $\bbeta$ is deterministic and satisfies
\[
    \frac{1}{|I_b|}\|\bbeta^{(b)}\|^2
    \rightarrow v_b\in[0,\infty),
    \qquad b=1,\ldots,k.
\]
Fix $\kappa>0$ and a positive integer $q$.
For any deterministic choice of distinct coordinates
$j_1,\ldots,j_q$, with $j_r\in I_{a_r}$ and fixed
$a_1,\ldots,a_q$, assume additionally that
\[
    \sup_N|\beta_{j_r}|<\infty,
    \qquad r=1,\ldots,q.
\]
Then
\begin{equation}\label{eq:joint_clt}
    \left(
        \mathcal C_{a_r}\hat\beta_{j_r}-\beta_{j_r}
    \right)_{r=1}^{q}
    \xrightarrow{\mathrm d}
    \mathcal N_q\!\left(
        0,
        \operatorname{diag}\bigl(
            V_{a_1}(\kappa),\ldots,V_{a_q}(\kappa)
        \bigr)
    \right),
\end{equation}
where
\begin{equation}\label{eq:coordinate_variance}
\begin{aligned}
V_a(\kappa)
:=\lim_{N,p\to\infty}
\mathcal C_a^2\kappa^2m_a(\kappa)^2
\Bigg[
&\sum_{b=1}^k
v_b\big[(\bI_k-\Omega(\kappa))^{-1}\Omega(\kappa)\big]_{ab}
\\
&+\sigma^2\sum_{b=1}^k
\big[(\bI_k-\Omega(\kappa))^{-1}\big]_{ab}
\frac1p\Tr\!\left(
\tilde\bT_p(\kappa)\Lambda_b\tilde\bT_p(\kappa)
\right)
\Bigg].
\end{aligned}
\end{equation}
\end{theorem}
\begin{remark}
The joint CLT is not restricted to the AR(1) covariance structure: it extends to deterministic positive semidefinite matrices $\Lambda_a$ with uniformly bounded operator norms, provided the remaining assumptions hold, the debiasing factors remain bounded, and the limits defining $V_a(\kappa)$ exist.
For the particular model \eqref{diagonal_lambda_model}, with covariance matrices given by \eqref{eq:def_Lambda}, these limits can be evaluated explicitly; see
Proposition~\ref{prop:explicit_variance_limit}.
\end{remark}

\begin{proof}[Sketch of the proof of Theorem \ref{thm:joint_clt}]
Here we introduce the main idea of the proof, the full proof can be found in Appendix~\ref{apx:proof_join_clt}.

Let $J=\{j_1,\ldots,j_q\}$ be the set of coordinates in the statement
of the theorem, with $j_r\in I_{a_r}$. We write
$\bbeta_J=(\beta_{j_1},\ldots,\beta_{j_q})^\intercal$ and
$\widehat{\bbeta}_J=(\hat\beta_{j_1},\ldots,\hat\beta_{j_q})^\intercal$,
and let $\mathbf C_J=\operatorname{diag}(\mathcal C_{a_1},\ldots,
\mathcal C_{a_q})$. Our aim is to show that the vector of debiased errors
$\mathbf C_J\widehat{\bbeta}_J-\bbeta_J$ converges to a centered Gaussian
vector with diagonal covariance matrix and with $r$th variance
$V_{a_r}(\kappa)$. To do this, we first show that
\begin{align}\label{eq:condit_var}
     \mathbf C_J\widehat{\bbeta}_J-\bbeta_J
    =\mathbf G_{J,p}+o_p(1),
\end{align}
where, conditionally on the design columns outside $J$ and on $\Z_N$,
the leading term $\mathbf G_{J,p}$ is centered Gaussian with covariance
matrix
\[
    \operatorname{diag}\bigl(
        V_{a_1,p}^{(-J)}(\kappa),\ldots,V_{a_q,p}^{(-J)}(\kappa)
    \bigr).
\]
It then remains to show that, for each $r=1,\ldots,q$,
\[
    V_{a_r,p}^{(-J)}(\kappa)
    \xrightarrow{\mathbb P} V_{a_r}(\kappa),
\]
where $V_{a_r}(\kappa)$ is the deterministic variance given in
\eqref{eq:coordinate_variance}.

In order to obtain \eqref{eq:condit_var}, we use the classical idea of leave-one-out. Let $\pphi_l$ denote the $l$th column of $p^{-1/2}\PPhi$ and write
$\PPhi^{(J)}=(\pphi_{j_1},\ldots,\pphi_{j_q})$. We separate the columns
indexed by $J$ from the rest by defining
\[
    \tilde\bQ_p^{(-J)}(\kappa)
    :=\left(\sum_{l\notin J}\pphi_l\pphi_l^\intercal
             +\kappa\bI_N\right)^{-1},
    \qquad
    \Y_N^{(-J)}:=\sum_{l\notin J}\pphi_l\beta_l+\Z_N.
\]
Both quantities are independent of the selected columns. By Woodbury's identity,
\[
    \widehat{\bbeta}_J-\bbeta_J
    =\kappa[\bQ_p(\kappa)]_{J,J}
    \left(
        \PPhi^{(J)\intercal}\tilde\bQ_p^{(-J)}(\kappa)\Y_N^{(-J)}
        -\bbeta_J
    \right).
\]
Using the deterministic equivalent of the resolvent \eqref{eq:gram_det_eq}, for $r,s=1,\ldots,q$, we can replace $[\bQ_p(\kappa)]_{j_r,j_s}=\delta_{rs}m_{a_r}(\kappa)+o_p(1).$
Since $q$ is fixed and the vector in parentheses is $O_p(1)$,
we obtain, jointly for $r=1,\ldots,q$,
\[
    \hat\beta_{j_r}
    =(1-\kappa m_{a_r}(\kappa))\beta_{j_r}
    +\kappa m_{a_r}(\kappa)\,
     \pphi_{j_r}^{\intercal}\tilde\bQ_p^{(-J)}(\kappa)\Y_N^{(-J)}
    +o_p(1).
\]
Multiplying by $\mathcal C_{a_r}=(1-\kappa m_{a_r}(\kappa))^{-1}$
and subtracting $\beta_{j_r}$ gives
\[
    \mathcal C_{a_r}\hat\beta_{j_r}-\beta_{j_r}
    =\mathcal C_{a_r}\kappa m_{a_r}(\kappa)\,
     \pphi_{j_r}^{\intercal}\tilde\bQ_p^{(-J)}(\kappa)\Y_N^{(-J)}
    +o_p(1).
\]

The leading terms on the right-hand side form the vector $\mathbf G_{J,p}$
in \eqref{eq:condit_var}. Conditionally on the columns outside $J$ and
on $\Z_N$, the vector $\tilde\bQ_p^{(-J)}(\kappa)\Y_N^{(-J)}$ is fixed,
while the selected columns $\pphi_{j_1},\ldots,\pphi_{j_q}$ remain
independent, with $\pphi_{j_r}\sim\mathcal N(0,\Lambda_{a_r}/p)$.
Thus, the coordinates of $\mathbf G_{J,p}$ are conditionally independent
and centered Gaussian. For $j_r\in I_a$, the corresponding conditional
variance is
\[
    V_{a,p}^{(-J)}(\kappa)
    =\frac{\mathcal C_a^2\kappa^2m_a(\kappa)^2}{p}
      (\Y_N^{(-J)})^\intercal
      \tilde\bQ_p^{(-J)}(\kappa)\Lambda_a
      \tilde\bQ_p^{(-J)}(\kappa)\Y_N^{(-J)}.
\]

To identify the limit of this variance, we first prove that is concentrates around its expectation by using the Gaussian Poincar\'e inequality.
In order to calculate the expectation we expand 
$\Y_N^{(-J)}=\sum_{l\notin J}\pphi_l\beta_l+\Z_N$, and notice that the mixed terms with noise and signal have zero expectation.
For the noise term, taking the expectation with respect to $\Z_N$
gives a trace. For the signal term, we remove one additional column
so that we can take the expectation with respect to this column
independently of the remaining resolvent. Applying the double-resolvent
deterministic equivalent in Lemma~\ref{lem:double_resolvent} then gives
the two contributions to $V_a(\kappa)$ in \eqref{eq:coordinate_variance}.
Hence $V_{a,p}^{(-J)}(\kappa)\xrightarrow{\mathbb P}V_a(\kappa)$.

Finally, for every $\mathbf t=(t_1,\ldots,t_q)^\intercal\in\mathbb R^q$,
the conditional Gaussian representation gives
\[
    \mathbb E\exp\!\left(i\mathbf t^\intercal\mathbf G_{J,p}\right)
    =\mathbb E\exp\!\left(
        -\frac12\sum_{r=1}^q t_r^2 V_{a_r,p}^{(-J)}(\kappa)
      \right)
    \longrightarrow
    \exp\!\left(-\frac12\sum_{r=1}^q t_r^2V_{a_r}(\kappa)\right).
\]
This proves the joint Gaussian limit for $\mathbf G_{J,p}$, and
\eqref{eq:condit_var} gives the same limit for the debiased errors
by Slutsky's theorem.
\end{proof}

\section{Consistent estimation of the limiting variances}
\label{sec:consistent_variance}

The limiting variances in Theorem~\ref{thm:joint_clt} depend on the
blockwise signal energies $v_1,\ldots,v_k$ and the noise variance
$\sigma^2$. We employ the method of moments (mom): we estimate these quantities jointly by Gaussian moment
equations and substitute the estimates directly into
\eqref{eq:coordinate_variance}. In particular, neither the construction
of the energy estimates nor that of the noise-variance estimate requires
ridge-based moment equations or a preliminary estimate of the other
unknown quantities.

Throughout, $k$ is fixed, $p_b:=|I_b|$,
$\ell_{b,p}:=p_b/p\to\ell_b>0$, and $c_p:=N/p\to c\in(0,\infty)$.
As in Theorem~\ref{thm:joint_clt}, the signal is deterministic and
\[
 v_{b,p}:=\frac{\|\bbeta^{(b)}\|^2}{p_b}\rightarrow v_b\in[0,\infty),
 \qquad \boldsymbol v_p:=(v_{1,p},\ldots,v_{k,p})^\top.
\]
We assume that the known covariance matrices are positive semidefinite
and satisfy, for constants $M<\infty$ and $\mu_->0$,
\begin{equation}\label{eq:var_est_trace_lower}
 \max_b\|\Lambda_b\|\le M,
 \qquad
 \min_b\frac1N\Tr\Lambda_b\ge\mu_-
\end{equation}
for all sufficiently large $p$. These conditions will be verified below
for the stationary AR(1) covariance matrices in \eqref{eq:def_Lambda}.

\subsection{Joint moment equations.}
For each block $b$, define $\PPhi_{I_b}=\LLambda_{b,N}^{1/2}\mathbf{G}_N^{(b)}$ (see \eqref{diagonal_lambda_model}) and 
\begin{equation}\label{eq:var_est_observed_moments}
 \bS_{b,p}:=\frac1p\PPhi_{I_b}\PPhi_{I_b}^\top,
 \qquad
 u_p:=\frac{\|\Y_N\|^2}{N},
 \qquad
 w_{b,p}:=\frac1N\Y_N^\top\bS_{b,p}\Y_N.
\end{equation}
Introduce the deterministic trace quantities
\begin{equation}\label{eq:var_est_moment_traces}
 \mu_b:=\frac1N\Tr\Lambda_b,
 \qquad
 \omega_{ba}:=\frac1N\Tr(\Lambda_b\Lambda_a),
 \qquad
 L_p:=\operatorname{diag}(\ell_{1,p},\ldots,\ell_{k,p}),
\end{equation}
and set
\begin{align}
 r_{b,p}^{\rm mom}
 &:=\frac{w_{b,p}}{\ell_{b,p}}-\mu_bu_p,
 &\boldsymbol r_p^{\rm mom}
 &:=(r_{1,p}^{\rm mom},\ldots,r_{k,p}^{\rm mom})^\top,
 \label{eq:var_est_noise_moments}\\
 \mathsf D_{p,ba}^{\rm mom}
 &:=\ell_{a,p}(\omega_{ba}-\mu_b\mu_a)
  +\mathbf1_{\{a=b\}}
       \left(c_p\mu_b^2+\frac{\omega_{bb}}p\right),
 &\mathsf D_p^{\rm mom}
 &:=(\mathsf D_{p,ba}^{\rm mom})_{b,a=1}^k.
 \label{eq:var_est_mom_matrix}
\end{align}

\begin{lemma}[Exact moments and concentration]
\label{lem:var_est_noise_moments}
Under the preceding assumptions,
\begin{align}
 \mathbb E[u_p]
 &=\sigma^2+\sum_{a=1}^k\ell_{a,p}\mu_av_{a,p},
 \label{eq:var_est_Eu}\\
 \mathbb E[\boldsymbol r_p^{\rm mom}]
 &=\mathsf D_p^{\rm mom}\boldsymbol v_p.
 \label{eq:var_est_Er}
\end{align}
Moreover,
\begin{equation}\label{eq:var_est_moment_concentration}
 \operatorname{Var}(u_p)=O(p^{-1}),
 \qquad
 \max_{1\le b\le k}\operatorname{Var}(w_{b,p})=O(p^{-1}),
 \qquad
 \max_{1\le b\le k}\operatorname{Var}(r_{b,p}^{\rm mom})=O(p^{-1}).
\end{equation}
\end{lemma}
\begin{proof}
Write the unscaled columns of $\PPhi$ as
$z_j\sim\mathcal N(0,\Lambda_{a(j)})$, independently, where
$a(j)=b$ for $j\in I_b$. Then
\[
 \Y_N=\frac1{\sqrt p}\sum_{j=1}^p\beta_jz_j+\Z_N,
 \qquad
 \Sigma_{Y,p}:=\operatorname{Cov}(\Y_N)
  =\sigma^2\bI_N+\sum_{a=1}^k\ell_{a,p}v_{a,p}\Lambda_a.
\]
Thus $\Y_N$ is centered Gaussian and $\|\Sigma_{Y,p}\|=O(1)$.
Taking the trace proves \eqref{eq:var_est_Eu}, and the Gaussian
quadratic-form identity gives
\[
 \operatorname{Var}(u_p)
 =\frac{2}{N^2}\Tr(\Sigma_{Y,p}^2)=O(N^{-1})=O(p^{-1}).
\]

For $i\in I_b$, the pair $(z_i,\Y_N)$ is jointly centered Gaussian,
with
$\operatorname{Cov}(z_i,\Y_N)=\beta_i\Lambda_b/\sqrt p$.
Therefore we obtain,
\begin{equation}\label{eq:var_est_pairing}
 \mathbb E(z_i^\top\Y_N)^2
 =\Tr(\Lambda_b\Sigma_{Y,p})
 +\frac{\beta_i^2}{p}
       \left\{(\Tr\Lambda_b)^2+\Tr(\Lambda_b^2)\right\}.
\end{equation}
Indeed, 
expanding $\sum_{r,s}\mathbb E[z_{i,r}Y_rz_{i,s}Y_s]$
and using Wick's formula yields the above result.
Since
$w_{b,p}=(Np)^{-1}\sum_{i\in I_b}(z_i^\top\Y_N)^2$,
summing \eqref{eq:var_est_pairing} yields the exact identity
\begin{equation}\label{eq:var_est_Ew}
 \mathbb E[w_{b,p}]
 =\ell_{b,p}\left[
   \sigma^2\mu_b+\sum_{a=1}^k\ell_{a,p}\omega_{ba}v_{a,p}
   +\left(c_p\mu_b^2+\frac{\omega_{bb}}p\right)v_{b,p}
  \right].
\end{equation}
Subtracting $\mu_b\mathbb E[u_p]$ after division by $\ell_{b,p}$
proves \eqref{eq:var_est_Er}.

To prove concentration of $w_{b,p}$, write
$z_j=\Lambda_{a(j)}^{1/2}g_j$ and $\Z_N=\sigma g_0$ with all
$g_j,g_0$ independent standard Gaussian vectors. The Gaussian
Poincar\'e inequality gives
\[
 \operatorname{Var}(w_{b,p})
 \le\sum_{j=0}^p\mathbb E\|\nabla_{g_j}w_{b,p}\|^2.
\]
Direct differentiation gives
\begin{align*}
 \nabla_{g_j}w_{b,p}
 &=\frac{2\beta_j}{N\sqrt p}\Lambda_{a(j)}^{1/2}\bS_{b,p}\Y_N
  +\frac{2\mathbf1_{\{j\in I_b\}}}{Np}
        \Lambda_b^{1/2}\Y_N(z_j^\top\Y_N),\\
 \nabla_{g_0}w_{b,p}&=\frac{2\sigma}{N}\bS_{b,p}\Y_N.
\end{align*}
Using $\|x+y\|^2\le2\|x\|^2+2\|y\|^2$ and
\[
 \sum_{j\in I_b}(z_j^\top\Y_N)^2
 =p\Y_N^\top\bS_{b,p}\Y_N
 \le p\|\bS_{b,p}\|\|\Y_N\|^2,
\]
we obtain
\begin{align}\label{eq:var_est_gradient_bound}
 \operatorname{Var}(w_{b,p})
 \le C\bigg[&\frac{\|\bbeta\|^2/p+\sigma^2}{N^2}
       \mathbb E\bigl(\|\bS_{b,p}\|^2\|\Y_N\|^2\bigr)
       \nonumber\\
 &+\frac1{N^2p}
       \mathbb E\bigl(\|\bS_{b,p}\|\|\Y_N\|^4\bigr)\bigg].
\end{align}
Note that $\|\Sigma_{Y,p}\|=O(1)$,
$\mathbb E\|\Y_N\|^4=O(N^2)$ and
$\mathbb E\|\Y_N\|^8=O(N^4)$.
Further, we have
\[
 \|\bS_{b,p}\|\le\frac{M}{p}\|G^{(b)}\|^2,
 \qquad
 \mathbb E\|G^{(b)}\|^8=O\bigl((N+p_b)^4\bigr),
\]
so $\sup_p\mathbb E\|\bS_{b,p}\|^4<\infty$.

Cauchy-Schwarz now bounds the two expectations in
\eqref{eq:var_est_gradient_bound} by $O(N)$ and $O(N^2)$,
respectively, without assuming independence of $\bS_{b,p}$ and
$\Y_N$. Since $N\asymp p$, the desired $O(p^{-1})$ bound follows.
Finally, $\ell_{b,p}$ is bounded away from zero and $\mu_b\le M$;
applying the variance inequality to
$r_{b,p}^{\rm mom}=w_{b,p}/\ell_{b,p}-\mu_bu_p$ proves the last
claim in \eqref{eq:var_est_moment_concentration}.
\end{proof}

\subsection{Identification and simultaneous estimation.}
Define the raw estimators
\begin{equation}\label{eq:var_est_joint_raw}
 \widetilde{\boldsymbol v}_p
 :=(\mathsf D_p^{\rm mom})^{-1}\boldsymbol r_p^{\rm mom},
 \qquad
 \widetilde\sigma_p^2
 :=u_p-\sum_{a=1}^k\ell_{a,p}\mu_a\widetilde v_{a,p}.
\end{equation}
Their nonnegative versions are
\begin{equation}\label{eq:var_est_joint_estimators}
 \widehat v_b:=[\widetilde v_{b,p}]_+,
 \qquad
 \widehat\sigma^2:=[\widetilde\sigma_p^2]_+,
 \qquad [x]_+:=\max\{x,0\},
\end{equation}
and
$\widehat{\boldsymbol v}:=(\widehat v_1,\ldots,\widehat v_k)^\top$.
Thus the noise estimate in \eqref{eq:var_est_joint_estimators} is
computed using the raw energy estimates before truncation.

\begin{proposition}[Uniform identification and consistency]
\label{prop:var_est_joint_moments}
There exists $p_0$, such that under \eqref{eq:var_est_trace_lower}, the moment matrix is invertible
and satisfies for 
\begin{equation}\label{eq:var_est_inverse_bound}
 \sup_{p \geq p_0}\|(\mathsf D_p^{\rm mom})^{-1}\|<\infty.
\end{equation}
The raw estimators are exactly unbiased for $\boldsymbol v_p$ and
$\sigma^2$. Moreover,
\begin{equation}\label{eq:var_est_joint_MSE}
 \mathbb E\|\widetilde{\boldsymbol v}_p-\boldsymbol v_p\|^2
 +\mathbb E|\widetilde\sigma_p^2-\sigma^2|^2=O(p^{-1}),
\end{equation}
and the same bound holds for their nonnegative versions.
In particular,
\begin{equation}\label{eq:var_est_joint_consistency}
 \max_b|\widehat v_b-v_{b,p}|=O_p(p^{-1/2}),
 \qquad
 \widehat\sigma^2-\sigma^2=O_p(p^{-1/2}),
 \qquad
 \widehat v_b\xrightarrow{\mathbb P}v_b.
\end{equation}
\end{proposition}
\begin{proof}
Let $\mathsf K_p=(\omega_{ba}-\mu_b\mu_a)_{b,a=1}^k$. For every
$x\in\mathbb R^k$,
\[
 x^\top\mathsf K_px
 =\frac1N\Tr\left[
       \left\{\sum_{b=1}^kx_b(\Lambda_b-\mu_b\bI_N)\right\}^2
      \right]\ge0.
\]
Thus $\mathsf K_p$ is a symmetric positive semidefinite Gram matrix.
Since
\[
 \mathsf D_p^{\rm mom}
 =\mathsf K_pL_p+
    \operatorname{diag}\left(c_p\mu_b^2+\frac{\omega_{bb}}p\right)_{b=1}^k,
\]
its diagonal similarity transform satisfies
\begin{equation}\label{eq:var_est_positive_matrix}
\begin{aligned}
 L_p^{1/2}\mathsf D_p^{\rm mom}L_p^{-1/2}
 &=L_p^{1/2}\mathsf K_pL_p^{1/2}
   +\operatorname{diag}\left(c_p\mu_b^2+\frac{\omega_{bb}}p\right)_{b=1}^k
   \\
 &\succeq c_p\mu_-^2\bI_k.
\end{aligned}
\end{equation}
Choose $c_->0$ such that $c_p\ge c_-$ for large $p$.
Writing $\ell_{\min,p}=\min_b\ell_{b,p}$ and
$\ell_{\max,p}=\max_b\ell_{b,p}$ gives the explicit bound
\[
 \|(\mathsf D_p^{\rm mom})^{-1}\|
 \le\frac1{c_-\mu_-^2}
       \sqrt{\frac{\ell_{\max,p}}{\ell_{\min,p}}}=O(1).
\]
This proves stable identification, even when $\mathsf K_p$ is singular.

The exact moment equations imply
$\mathbb E\widetilde{\boldsymbol v}_p=\boldsymbol v_p$ and
$\mathbb E\widetilde\sigma_p^2=\sigma^2$. Further,
\[
 \widetilde{\boldsymbol v}_p-\boldsymbol v_p
 =(\mathsf D_p^{\rm mom})^{-1}
       (\boldsymbol r_p^{\rm mom}-\mathbb E\boldsymbol r_p^{\rm mom}).
\]
The inverse bound and Lemma~\ref{lem:var_est_noise_moments} yield
$\mathbb E\|\widetilde{\boldsymbol v}_p-\boldsymbol v_p\|^2=O(p^{-1})$.
For $d_p=(\ell_{a,p}\mu_a)_{a=1}^k$, which has bounded norm,
\[
 \widetilde\sigma_p^2-\sigma^2
 =(u_p-\mathbb E u_p)
   -d_p^\top(\widetilde{\boldsymbol v}_p-\boldsymbol v_p).
\]
Together with concentration of $u_p$ by Lemma~\ref{lem:var_est_noise_moments} this proves \eqref{eq:var_est_joint_MSE}.  Since $|[x]_+-t|\le|x-t|$ for $t\ge0$,
truncation preserves the mean-square bounds. Chebyshev's inequality
and $v_{b,p}\to v_b$ complete the proof.
\end{proof}

\begin{remark}[Verification for the covariance model]
\label{rem:var_est_AR_identification}
Since the number of blocks and their autoregressive coefficients are
fixed, there is $r<1$ with $\max_b|\lambda_b|\le r$.
The absolute row-sum bound gives
\[
 \|\Lambda_b\|\le(1-|\lambda_b|)^{-2}\le(1-r)^{-2},
 \qquad \mu_b=(1-\lambda_b^2)^{-1}\ge1.
\]
Thus \eqref{eq:var_est_trace_lower} holds, and the lower bound in
\eqref{eq:var_est_positive_matrix} is at least $c_p\bI_k$.
In particular, the required identifiability condition holds in the
model of the paper.
For computation, the exact trace formula is
\[
 \omega_{ba}
 =\mu_b\mu_a\left[1+2\sum_{h=1}^{N-1}
                   \left(1-\frac hN\right)(\lambda_b\lambda_a)^h\right].
\]
\end{remark}

\subsection{Plug-in estimation of the CLT variances.}
Fix $\kappa>0$. Retain the finite-dimensional fixed-point solutions
$m_a(\kappa)$ and $\mathcal C_a=(1-\kappa m_a)^{-1}$ used in the
preceding CLT. To keep the estimator fully finite-dimensional, set
\begin{align}
 U_p&:=\bI_N+\sum_{b=1}^k\ell_{b,p}m_b\Lambda_b,
 &\tilde\bT_p&=\kappa^{-1}U_p^{-1},\nonumber\\
 (\Omega_p)_{ab}
 &:=\ell_{b,p}\kappa^2m_b^2\frac1p\Tr
           (\Lambda_a\tilde\bT_p\Lambda_b\tilde\bT_p),
 & B_p&:=(\bI_k-\Omega_p)^{-1},\label{eq:var_est_Omega_finite}\\
 \zeta_{b,p}&:=\frac1p\Tr(\tilde\bT_p\Lambda_b\tilde\bT_p),
 &\boldsymbol\zeta_p&:=(\zeta_{1,p},\ldots,\zeta_{k,p})^\top.
 \label{eq:var_est_zeta}
\end{align}
Dependence on $\kappa$ is suppressed on the right-hand sides.
The sole difference from the notation in \eqref{eq:coordinate_variance}
is the use of the observed block proportion $\ell_{b,p}$ rather than
its limit $\ell_b$ in $\Omega_p$.

\begin{lemma}[Uniform stability of the variance coefficients]
\label{lem:var_est_coefficient_stability}
For each fixed $\kappa>0$, the matrices $B_p$ and $B_p\Omega_p$,
the vector $\boldsymbol\zeta_p$, and the constants $\mathcal C_a$
are uniformly bounded for large $p$. Moreover, $\Omega_p$, $B_p$,
and $B_p\Omega_p$ are entrywise nonnegative, and there are constants
$m_*>0$ and $\zeta_*>0$ such that $m_a\ge m_*$ and
$\zeta_{a,p}\ge\zeta_*$ for every block and large $p$.
\end{lemma}
\begin{proof}
Choose $0<c_-\le c_p\le c_+<\infty$. Put
\[
 \delta_a:=\frac1p\Tr(\Lambda_aU_p^{-1}),
 \qquad s_a:=m_a^{-1}=\kappa+\delta_a,
 \qquad \boldsymbol s:=(s_a)_{a=1}^k.
\]
The fixed-point equations and $\|\Lambda_a\|\le M$ imply
\[
 \kappa\le s_a\le W_\kappa:=\kappa+c_+M,
 \qquad
 \bI_N\preceq U_p\preceq(1+M/\kappa)\bI_N.
\]
Consequently,
\[
 \delta_a\ge\delta_*:=\frac{c_-\mu_-}{1+M/\kappa}>0,
 \qquad
 m_a\ge W_\kappa^{-1},
 \qquad
 1\le\mathcal C_a=\frac{\kappa+\delta_a}{\delta_a}
                    \le\frac{W_\kappa}{\delta_*}.
\]
The entries of $\Omega_p$ are nonnegative, since they are positive
multiples of traces of products of positive semidefinite matrices
after conjugation by $U_p^{-1/2}$. Also,
\begin{align*}
 (\Omega_p\boldsymbol s)_a
 &=\frac1p\Tr\{\Lambda_aU_p^{-1}(U_p-\bI_N)U_p^{-1}\}\\
 &=\delta_a-\frac1p\Tr(\Lambda_aU_p^{-2})
 \le s_a-\kappa.
\end{align*}
For $W_p:=\operatorname{diag}(s_1,\ldots,s_k)$, it follows that
\[
 \|W_p^{-1}\Omega_pW_p\|_\infty
 =\max_a\frac{(\Omega_p\boldsymbol s)_a}{s_a}
 \le1-\frac{\kappa}{W_\kappa}<1.
\]
Here $\|\cdot\|_\infty$ is the maximum absolute row-sum norm.
Thus the spectral radius of $\Omega_p$ is uniformly less than one.
The Neumann series $B_p=\sum_{j\ge0}\Omega_p^j$, together with
the uniform bounds on $W_p$ and $W_p^{-1}$ and fixed $k$, proves
uniform boundedness and entrywise nonnegativity of $B_p$ and
$B_p\Omega_p$. Finally,
\[
 \frac{c_-\mu_-}{\kappa^2(1+M/\kappa)^2}
 \le\zeta_{a,p}
 =\frac1{p\kappa^2}\Tr(\Lambda_aU_p^{-2})
 \le\frac{c_+M}{\kappa^2},
\]
which proves the remaining assertions.
\end{proof}

Using the joint estimators in \eqref{eq:var_est_joint_estimators}, define
\begin{equation}\label{eq:var_est_Vhat}
 \widehat V_a(\kappa)
 :=\mathcal C_a^2\kappa^2m_a(\kappa)^2
   \left\{
       [B_p\Omega_p\widehat{\boldsymbol v}]_a
       +\widehat\sigma^2[B_p\boldsymbol\zeta_p]_a
   \right\},
 \qquad a=1,\ldots,k.
\end{equation}
All coefficients in this formula are computed from the known covariance
matrices, the finite-dimensional fixed-point system, and the observed
block sizes. The only statistically estimated quantities are
$\widehat{\boldsymbol v}$ and $\widehat\sigma^2$. By the nonnegativity
in Lemma~\ref{lem:var_est_coefficient_stability},
$\widehat V_a(\kappa)\ge0$ for every sample.

\begin{theorem}[Variance consistency and feasible joint inference]
\label{thm:var_est_studentized}
Assume the preceding moment-model conditions and the hypotheses of
Theorem~\ref{thm:joint_clt}, including existence of the limits
$V_a(\kappa)$ in \eqref{eq:coordinate_variance}. For each fixed
$\kappa>0$,
\begin{equation}\label{eq:var_est_V_consistency}
 \max_{1\le a\le k}|\widehat V_a(\kappa)-V_a(\kappa)|
 \xrightarrow{\mathbb P}0.
\end{equation}
If $\sigma^2>0$, all these limiting variances are strictly positive.
More generally, whenever $V_{a_r}(\kappa)>0$ for the fixed
coordinates $j_r\in I_{a_r}$ under consideration, define
$\widehat V_{a,+}:=\widehat V_a\vee\varepsilon_p$ for any deterministic
$\varepsilon_p>0$ with $\varepsilon_p\downarrow0$. Then
\begin{equation}\label{eq:var_est_joint_studentization}
 \left(
   \frac{\mathcal C_{a_r}\widehat\beta_{j_r}-\beta_{j_r}}
        {\sqrt{\widehat V_{a_r,+}(\kappa)}}
 \right)_{r=1}^q
 \xrightarrow{\mathrm d}\mathcal N_q(0,\bI_q).
\end{equation}
Consequently, for any one of these selected coordinates $j\in I_a$ and $0<\alpha<1$, the interval
\begin{equation}\label{eq:var_est_confidence_interval}
 \left[
  \mathcal C_a\widehat\beta_j
       -z_{1-\alpha/2}\sqrt{\widehat V_{a,+}(\kappa)},\quad
  \mathcal C_a\widehat\beta_j
       +z_{1-\alpha/2}\sqrt{\widehat V_{a,+}(\kappa)}
 \right]
\end{equation}
has coverage tending to $1-\alpha$. Moreover,
\begin{equation}\label{eq:var_est_Wald}
 \sum_{r=1}^q
   \frac{(\mathcal C_{a_r}\widehat\beta_{j_r}-\beta_{j_r})^2}
        {\widehat V_{a_r,+}(\kappa)}
 \xrightarrow{\mathrm d}\chi_q^2.
\end{equation}
Here $z_{1-\alpha/2}$ is the corresponding standard-normal quantile.
\end{theorem}
\begin{proof}
Define the finite-dimensional deterministic variance expression
\[
 V_{a,p}^{\rm det}(\kappa)
 :=\mathcal C_a^2\kappa^2m_a^2
    \left\{[B_p\Omega_p\boldsymbol v_p]_a
               +\sigma^2[B_p\boldsymbol\zeta_p]_a\right\}.
\]
The coefficients are uniformly bounded by
Lemma~\ref{lem:var_est_coefficient_stability}. Therefore,
\begin{equation}\label{eq:var_est_plugin_error}
 \max_a|\widehat V_a-V_{a,p}^{\rm det}|
 \le C_\kappa\left(
       \|\widehat{\boldsymbol v}-\boldsymbol v_p\|
          +|\widehat\sigma^2-\sigma^2|\right)
 =O_p(p^{-1/2}).
\end{equation}

Replacing $\boldsymbol v_p$ by its limit $\boldsymbol v$ changes
$V_{a,p}^{\rm det}$ by $o(1)$. Likewise, replacing $\ell_{b,p}$ by
$\ell_b$ in $\Omega_p$ changes that matrix by $o(1)$, because its
remaining factors are bounded. The uniform inverse bound and the
resolvent identity for matrix inverses show that the corresponding
inverse matrices also differ by $o(1)$. Thus
$V_{a,p}^{\rm det}(\kappa)\to V_a(\kappa)$, exactly as defined
in \eqref{eq:coordinate_variance}. Combining this with
\eqref{eq:var_est_plugin_error} proves variance consistency.

For positivity, $B_p\boldsymbol\zeta_p\ge\boldsymbol\zeta_p$
componentwise, $\mathcal C_a\ge1$, and all signal terms are
nonnegative. Hence, when $\sigma^2>0$,
\[
 V_{a,p}^{\rm det}(\kappa)
 \ge\sigma^2\kappa^2m_*^2\zeta_*>0
\]
uniformly for large $p$. Passing to the limit proves $V_a(\kappa)>0$.

Finally, let
\[
 \widehat\Sigma_{J,p}
 :=\operatorname{diag}(\widehat V_{a_r,+})_{r=1}^q,
 \qquad
 \Sigma_J:=\operatorname{diag}(V_{a_r})_{r=1}^q.
\]
Variance consistency implies
$\widehat\Sigma_{J,p}^{-1/2}\to_{\mathbb P}\Sigma_J^{-1/2}$
whenever the displayed limiting variances are positive.
Theorem~\ref{thm:joint_clt} and multivariate Slutsky then give
\eqref{eq:var_est_joint_studentization}. The interval coverage follows from
continuity of the standard-normal distribution, and
\eqref{eq:var_est_Wald} follows by the continuous mapping theorem.
\end{proof}

\section{Connection with Bayesian estimation }
\label{sec:bayes_conn}

The deterministic-equivalent equations underlying ridge regression have
an exact interpretation in the Bayesian formulation of the same block
model. For a centered Gaussian signal prior and a suitable choice of the
ridge parameter, the ridge estimator becomes the posterior mean.
The blockwise resolvent then determines the Bayesian estimation errors,
and its limiting fixed point coincides with the Gaussian specialization
of the scalar-channel equations in \cite{mutinf_time}. We establish this
correspondence below, including its implications for the block-MMSE
formula and the variance of the debiased estimator.
For a Gaussian prior, we prove the blockwise MMSE formula conjectured in
\cite[Eq.~(II.15)]{mutinf_time}. For the interested reader, in Appendix~\ref{app:bayes_additional},
we use the same equations to recover the Gaussian formulas for the
measurement MMSE and the mutual information per signal coordinate.

Throughout this section, $k$ is fixed, $p_a:=|I_a|$,
$\ell_{a,p}:=p_a/p\to\ell_a>0$, and $N/p\to c\in(0,\infty)$.
We retain the covariance assumptions and deterministic equivalents of
the preceding sections. For the explicit integral formulas we use the
stationary AR(1) covariance model and
Proposition~\ref{prop:explicit_variance_limit}. We write
$m_{a,p}(\kappa)$ when the dimension dependence of $m_a(\kappa)$
needs to be explicit, and $\bar m_a(\kappa)$ for its limit.
Set
\[
 X_p:=p^{-1/2}\PPhi,\qquad S_p:=X_p^\top X_p,\qquad
 Q_p(\kappa):=(S_p+\kappa\bI_p)^{-1},\qquad
 \widehat{\bbeta}_\kappa:=Q_p(\kappa)X_p^\top\Y_N.
\]
Thus $\kappa$ is the diagonal shift in the resolvent, equivalently the
penalty parameter in
$\frac12\|\Y_N-X_pu\|^2+\frac\kappa2\|u\|^2$.
Let $P_a$ denote the orthogonal coordinate projector onto $I_a$.

For the Bayesian statements, we replace the deterministic signal by
\begin{equation}\label{eq:bayes_conn_prior}
 \bbeta\sim\mathcal N(0,\rho\bI_p),\qquad
 \Y_N=X_p\bbeta+\Z_N,\qquad
 \Z_N\sim\mathcal N(0,\sigma^2\bI_N),
 \qquad \rho>0,\quad\sigma^2>0,
\end{equation}
with the signal, noise, and design mutually independent.

\subsection{Gaussian conjugacy and the blockwise MMSE.}

We choose the regularization parameter so that the ridge estimator
coincides with the posterior mean. For the Gaussian prior with variance
$\rho$ and noise variance $\sigma^2$, this choice is
\begin{equation}\label{eq:bayes_conn_matched_lambda}
    \kappa_{\mathrm B}:=\frac{\sigma^2}{\rho}.
\end{equation}
With this choice, for every $N,p$, the posterior distribution is
\begin{equation}\label{eq:bayes_conn_posterior}
 \mathcal L(\bbeta\mid\Y_N,\PPhi)
 =\mathcal N\!\left(
    \widehat{\bbeta}_{\kappa_{\mathrm B}},
    \sigma^2 Q_p(\kappa_{\mathrm B})
   \right).
\end{equation}
Indeed, conditional on the observed design and response, the posterior
density at $u$ is proportional to
\[
 \exp\left\{
 -\frac{1}{2\sigma^2}
 \left(\|\Y_N-X_pu\|^2+\frac{\sigma^2}{\rho}\|u\|^2\right)
 \right\}.
\]
Thus, maximizing the posterior density is equivalent to minimizing the
ridge objective with penalty $\kappa_{\mathrm B}$. The posterior
precision matrix is
\[
 \sigma^{-2}S_p+\rho^{-1}\bI_p
 =\sigma^{-2}(S_p+\kappa_{\mathrm B}\bI_p).
\]
Completing the square gives \eqref{eq:bayes_conn_posterior}. In
particular, the ridge estimator $\widehat{\bbeta}_{\kappa_{\mathrm B}}$,
being the posterior mean, minimizes the squared-error Bayes risk.

Next we define the blockwise and total parameter MMSEs by
\begin{align}
 \operatorname{mmse}_{a,p}
 &:=\frac1{p_a}\mathbb E
   \left\|P_a\{\bbeta-\mathbb E(\bbeta\mid\Y_N,\PPhi)\}\right\|^2,
 \label{eq:bayes_conn_block_mmse_def}\\
 \operatorname{mmse}_{p}
 &:=\frac1p\mathbb E
   \left\|\bbeta-\mathbb E(\bbeta\mid\Y_N,\PPhi)\right\|^2
   =\sum_a\ell_{a,p}\operatorname{mmse}_{a,p}.
 \nonumber
\end{align}

\begin{proposition}[Limiting block errors]
\label{prop:bayes_conn_posterior}
The blockwise and total MSEs admit the following limits:
\begin{align}
 \operatorname{mmse}_{a,p}
 &=\frac{\sigma^2}{p_a}\mathbb E\Tr\{P_aQ_p(\kappa_{\mathrm B})\}
   \rightarrow e_a:=\sigma^2\bar m_a(\kappa_{\mathrm B}),
 \label{eq:bayes_conn_block_mmse}\\
 \operatorname{mmse}_{p}
 &\rightarrow\sum_a\ell_ae_a.
 \label{eq:bayes_conn_total_mmse}
\end{align}
The same block limits hold for the conditional Bayes risk given the
design, in probability and in $L^1$. Moreover,
\begin{equation}\label{eq:bayes_conn_realized_block_loss}
 \frac1{p_a}\left\|
 P_a(\widehat{\bbeta}_{\kappa_{\mathrm B}}-\bbeta)
 \right\|^2\rightarrow e_a
 \quad\hbox{in }L^2.
\end{equation}
\end{proposition}

\begin{proof}

The conditional expectation minimizes posterior squared error, so
\[
 \frac1{p_a}\mathbb E\!\left[
 \|P_a(\bbeta-\widehat{\bbeta}_{\kappa_{\mathrm B}})\|^2
 \mid\Y_N,\PPhi\right]
 =\frac{\sigma^2}{p_a}\Tr(P_aQ_p(\kappa_{\mathrm B})).
\]
The right-hand side does not depend on $\Y_N$ and is also the
conditional Bayes risk given $\PPhi$.
The blockwise trace deterministic equivalent gives
\[
 \frac1{p_a}\Tr(P_aQ_p(\kappa_{\mathrm B}))
 -m_{a,p}(\kappa_{\mathrm B})\xrightarrow{\mathbb P}0.
\]
The factor $p/p_a$ is bounded, so this follows from the normalized
trace equivalent already used in the risk calculation. Since
$0\preceq\sigma^2Q_p(\kappa_{\mathrm B})\preceq\rho\bI_p$,
the conditional risk belongs to $[0,\rho]$. Convergence in probability
therefore implies convergence in $L^1$ and $L^2$, proving
\eqref{eq:bayes_conn_block_mmse}--\eqref{eq:bayes_conn_total_mmse}.

For the realized loss, the centered posterior error is Gaussian.
Its conditional quadratic-form variance is
\[
 \operatorname{Var}\!\left(
 \left.\frac{\|P_a(\bbeta-\widehat{\bbeta}_{\kappa_{\mathrm B}})\|^2}{p_a}
 \right|\Y_N,\PPhi\right)
 =\frac{2\sigma^4}{p_a^2}
       \Tr\{(P_aQ_p(\kappa_{\mathrm B})P_a)^2\}
 \le\frac{2\rho^2}{p_a}.
\]
Combining this bound with $L^2$ convergence of the conditional mean
proves \eqref{eq:bayes_conn_realized_block_loss}.
\end{proof}

\subsection{Identification of the two fixed-point systems.}

Next, we identify our fixed point system with the one obtained in \cite{mutinf_time}. Define
\[
 f_a(\theta):=(1-2\lambda_a\cos\theta+\lambda_a^2)^{-1},
 \qquad
 \tau(\theta,\kappa)
 :=1+\sum_b\ell_b\bar m_b(\kappa)f_b(\theta).
\]
By evenness of the generating functions, the limiting system in
\eqref{eq:limiting_fixed_point} can be written as
\begin{equation}\label{eq:bayes_conn_ridge_fp}
 \frac1{\bar m_a(\kappa)}
 =\kappa+\frac c\pi\int_0^\pi
         \frac{f_a(\theta)}{\tau(\theta,\kappa)}\,d\theta.
\end{equation}
For a scalar prior $P_0$, let
\[
 \operatorname{mmse}_{P_0}(s)
 :=\mathbb E\left[
  \{B-\mathbb E(B\mid\sqrt{s}B+Z)\}^2
 \right],\qquad B\sim P_0,\quad Z\sim\mathcal N(0,1),
\]
where $B,Z$ are independent. The stationary equations in
\cite[Eq.~(II.5)]{mutinf_time}, with its generating functions
$\delta_a$ written as $f_a$, are
\begin{equation}\label{eq:bayes_conn_bayes_fp}
 \begin{split}
 r_{2,a}&=\operatorname{mmse}_{P_0}(r_{1,a}),\\
 r_{1,a}&=\frac c\pi\int_0^\pi
   \frac{f_a(\theta)}
        {\sigma^2+\sum_b\ell_bf_b(\theta)r_{2,b}}\,d\theta.
 \end{split}
\end{equation}
The parameters $r_{1,a}$ are scalar-channel precisions; the quantities
$r_{2,a}$ are the corresponding scalar estimation errors. 

\begin{proposition}[Gaussian fixed-point correspondence]
\label{prop:bayes_conn_fixed_point}
For $P_0=\mathcal N(0,\rho)$, the positive solutions of
\eqref{eq:bayes_conn_bayes_fp} and
\eqref{eq:bayes_conn_ridge_fp} at $\kappa=\kappa_{\mathrm B}$
are in one-to-one correspondence through
\begin{equation}\label{eq:bayes_conn_dictionary}
 \quad
 r_{2,a}=\sigma^2\bar m_a(\kappa_{\mathrm B}),\qquad
 r_{1,a}=\frac1{\sigma^2\bar m_a(\kappa_{\mathrm B})}
                  -\frac1\rho.
 \quad
\end{equation}
This solution is unique, $0<r_{2,a}<\rho$, and $r_{1,a}>0$.
Consequently,
\begin{equation}\label{eq:bayes_conn_conjecture_gaussian}
 \lim_{p\to\infty}\operatorname{mmse}_{a,p}=r_{2,a},
 \qquad
 \lim_{p\to\infty}\operatorname{mmse}_{p}
       =\sum_a\ell_ar_{2,a}.
\end{equation}
These are the Gaussian-prior versions of the block-MMSE formulas in
\cite[Eq.~(II.15)]{mutinf_time}.
\end{proposition}

\begin{proof}
The scalar Gaussian posterior variance is
\[
 \operatorname{mmse}_{\mathcal N(0,\rho)}(s)
 =\frac{\rho}{1+\rho s}=(\rho^{-1}+s)^{-1}.
\]
Thus the first Bayesian equation is
$r_{2,a}^{-1}=\rho^{-1}+r_{1,a}$.
Substituting $r_{2,a}=\sigma^2\bar m_a$ in the second equation gives
\[
 \sigma^2r_{1,a}
 =\frac c\pi\int_0^\pi
       \frac{f_a(\theta)}
            {1+\sum_b\ell_b\bar m_bf_b(\theta)}\,d\theta.
\]
Multiplying $r_{2,a}^{-1}=\rho^{-1}+r_{1,a}$ by $\sigma^2$
now gives \eqref{eq:bayes_conn_ridge_fp}. Conversely, any positive
solution of that ridge system defines $r_{1,a},r_{2,a}$ by
\eqref{eq:bayes_conn_dictionary} and satisfies both Bayesian equations.
The integral in \eqref{eq:bayes_conn_ridge_fp} is strictly positive,
so $\bar m_a<\kappa_{\mathrm B}^{-1}$, proving the asserted signs.

For completeness, uniqueness can be checked directly. Let
$F_a(u)$ denote the reciprocal of the right-hand side of
\eqref{eq:bayes_conn_ridge_fp}, with $\bar m$ replaced by $u$.
The map $F$ is coordinatewise increasing. For every $t>1$,
\[
 1+\sum_b\ell_btu_bf_b
 <t\left(1+\sum_b\ell_bu_bf_b\right)
 \quad\rightarrow\quad F_a(tu)<tF_a(u).
\]
If $u,v$ were distinct positive fixed points, after interchanging
them if necessary one could take $t=\max_a u_a/v_a>1$.
Then $u=F(u)\le F(tv)<tF(v)=tv$ in every coordinate,
contradicting equality in a coordinate attaining the maximum.
Existence follows from the resolvent fixed-point result in
Proposition~\ref{prop:explicit_variance_limit}. Finally,
Proposition~\ref{prop:bayes_conn_posterior} identifies each $r_{2,a}$
with the limiting block MMSE, not only with a scalar-channel parameter.
\end{proof}

\section{Numerical experiments}\label{sec:ridge_simulations}

In this section, we provide some numerical experiments to illustrate the
coordinate-wise normal approximation in Theorem~\ref{thm:var_est_studentized}.
We consider two choices of the block design in \eqref{diagonal_lambda_model},
both with $N/p=2$ and noise variance $\sigma^2=1$. For each design, we draw
the entries of $\bbeta$ independently from the Rademacher distribution and
keep this signal fixed over $2000$ realizations of the design matrix and
the noise.

For each realization, we compute the ridge estimator with $\kappa=5$ and
estimate the variance in each block using \eqref{eq:var_est_Vhat}. We take
the first coordinate $j$ in each block $I_a$ and compare the empirical
quantiles of
\[
    \frac{\mathcal C_a\widehat\beta_j-\beta_j}
         {\sqrt{\widehat V_a(\kappa)\vee p^{-1}}}
\]
with the quantiles of the standard normal distribution. We also plot the
theoretical variance as a function of $\kappa$. For these curves, we
evaluate the expression inside the limit in \eqref{eq:coordinate_variance}
at the dimensions used in the simulations, with $v_a=1$ for every block
and $\sigma^2=1$.

\begin{figure*}[!t]
    \centering
    \includegraphics[width=0.82\textwidth]{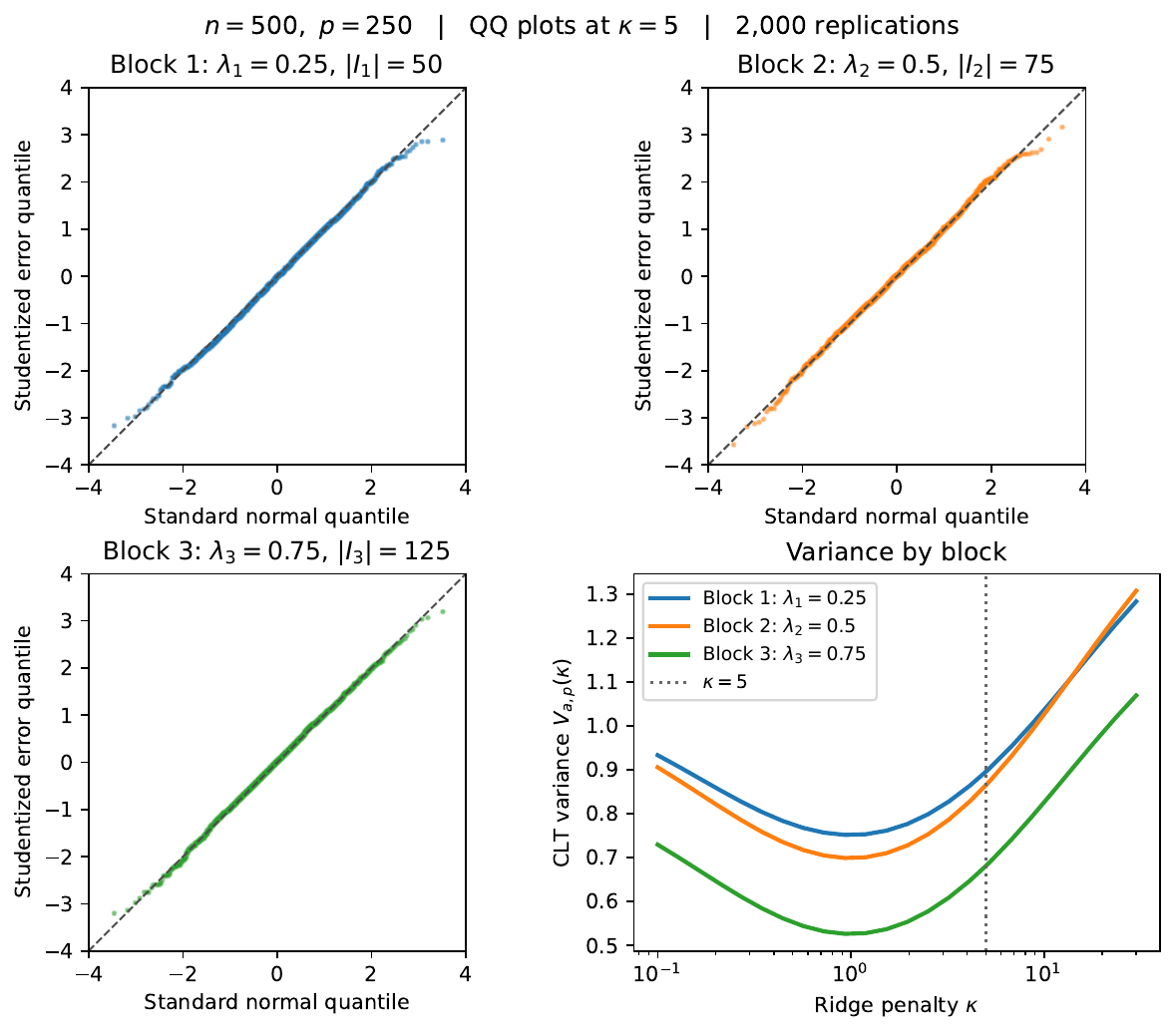}
    \caption{QQ plots for the first coordinate of each of the three blocks
    with $N=500$ and $p=250$. The upper panels correspond to the blocks
    of sizes $50$ and $75$, and the lower left panel to the block of size
    $125$. The lower right panel shows the theoretical variance for each
    block as a function of the ridge penalty. The dotted line marks
    $\kappa=5$, used in the QQ plots.}
    \label{fig:ridge_three_blocks}
\end{figure*}

Figure~\ref{fig:ridge_three_blocks} shows the results for $N=500$ and
$p=250$. We take three blocks of sizes $50$, $75$ and $125$, with
autoregressive parameters $1/4$, $1/2$ and $3/4$, respectively. In all
three QQ plots, the points follow the diagonal closely over most of the
range.  In the lower right panel, we see that
the third block has the smallest variance over the range of $\kappa$
shown. The curves for the first two blocks are closer and cross for
larger values of $\kappa$.

\begin{figure*}[!t]
    \centering
    \includegraphics[width=\textwidth]{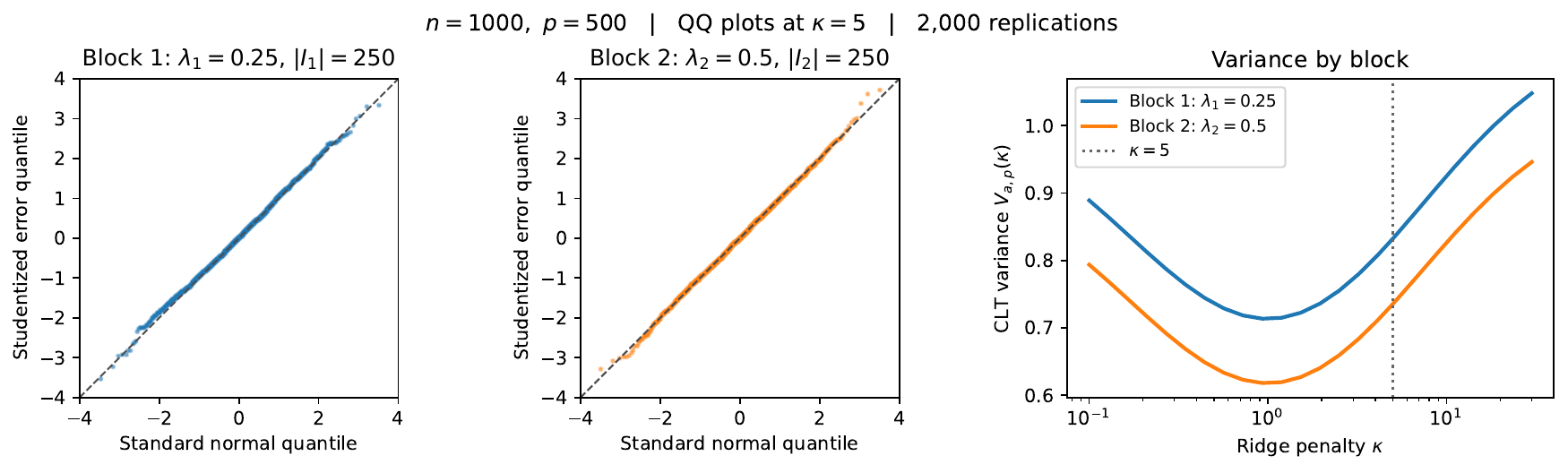}
    \caption{QQ plots and theoretical variances for $N=1000$, $p=500$
    and two blocks of size $250$. The left and middle panels correspond
    to autoregressive parameters $1/4$ and $1/2$, respectively. The right
    panel shows the two variance curves, with the dotted line marking
    the ridge penalty used in the QQ plots.}
    \label{fig:ridge_two_blocks}
\end{figure*}

In Figure~\ref{fig:ridge_two_blocks}, we take $N=1000$, $p=500$ and
two blocks of equal size, with autoregressive parameters $1/4$ and $1/2$.
For both blocks, the empirical quantiles again follow the diagonal closely,
with small deviations at the ends. The right panel shows that the block
with autoregressive parameter $1/2$ has the smaller variance throughout
the plotted range. In both figures, the variance curves reach their minima
near $\kappa=1$.

\section{Discussion}\label{sec:discussion}

In this paper, we develop debiased estimators of regression
coefficients under temporal dependence using tools from random matrix
theory. We characterize the asymptotic bias of the ridge estimator
through the theory of deterministic equivalents and use this
characterization to construct an asymptotically debiased estimator.
We establish its joint asymptotic normality, characterize the limiting
variances, and construct consistent estimators of these variances,
thereby enabling valid inference for the regression coefficients
without imposing sparsity assumptions on the underlying signal.

A central contribution of this work is the extension of debiasing
techniques from independent observations to temporally dependent
designs. Our framework provides a basis for studying other models
with temporal dependence, although the present analysis exploits
the particular structure of independent Gaussian design columns
with blockwise covariance matrices. This structure provides a
tractable starting point but restricts the class of models covered
by our results. An important direction is to accommodate observations
generated by more general dependent stochastic processes, allowing,
in particular, for dependence across covariates. Such an extension
could provide tools for inference in vector autoregressive models in the dense regime.

A second direction is to extend the analysis to non-Gaussian time
series. Without Gaussianity, blockwise right-rotational invariance
cannot generally be expected, even when the columns are independent.
More importantly, the conditional Gaussian representation that
underlies our proof of asymptotic normality is no longer available.
While universality results for random matrices may allow certain
aspects of the asymptotic risk characterization to be extended
beyond Gaussian designs, establishing the corresponding fluctuation
results presents a substantially different challenge. Addressing
non-Gaussian models will therefore require techniques that
accommodate temporal dependence while providing sufficiently
precise control of the fluctuations of the debiased estimator.

Finally, our results highlight a connection between Bayesian and
frequentist inference. Under a Gaussian prior and an appropriately
matched regularization parameter, the fixed-point equations obtained
through deterministic equivalents coincide with the Gaussian
specialization of those arising in the Bayesian formulation of
\cite{mutinf_time}. This correspondence not only establishes the
Gaussian specialization of the conjectured blockwise MMSE formula
but also provides an interpretation of the limiting variance of
the debiased estimator in terms of the effective scalar channel
in the Bayesian formulation. It would be interesting to determine
how broadly this correspondence holds, particularly for non-Gaussian
priors and other regularized estimators, and whether it can be used
to transfer principled inferential procedures between Bayesian
and frequentist settings.

\section*{Disclosure of AI Use}

GPT 6 was used to assist with reviewing the manuscript, checking mathematical arguments, and copyediting.
The authors have independently verified the content and take full responsibility for its correctness.

\bibliography{inf_hd_time}
\bibliographystyle{abbrv}

\appendix

\section{Deterministic equivalents}

The following definition of deterministic equivalent is adapted from \cite{Liao_RMT}

\begin{definition}\label{def_det_eq}
 We say that $\bar{Q}_n\in \mathbb{R}^{n\times n}$ is a deterministic equivalent for the symmetric random matrix $Q_n\in \mathbb{R}^{n\times n}$ if, for every sequence of deterministic matrices $A_n\in\mathbb{R}^{n\times n}$ with $\sup_n\|A_n\|<\infty$ and
deterministic unit vectors $a_n,b_n\in\mathbb{R}^n$, we have, as $n\to\infty$,

\begin{align*}
    \frac{1}{n}\tr[A_n(Q_n-\bar{Q}_n)]&\rightarrow 0\\
    a_n^T(Q_n-\bar{Q}_n)b_n&\rightarrow 0
\end{align*}
where the convergence is either in probability or almost sure.
\end{definition}
\begin{lemma}\label{lem:double_resolvent}
Let $\kappa>0$ be fixed. Define the $k\times k$ matrix
$\Omega(\kappa)$ by
\[
\Omega_{ab}(\kappa)
:=\ell_b\kappa^2m_b(\kappa)^2\frac1p\Tr\!\left(\Lambda_a\tilde\bT_p(\kappa)\Lambda_b\tilde\bT_p(\kappa)\right).
\]
Under our assumptions, $\bI_k-\Omega(\kappa)$ is invertible for sufficiently large $p$,
and
\[
\tilde\bQ_p(\kappa)\Lambda_a\tilde\bQ_p(\kappa)
\longleftrightarrow
\sum_{b=1}^k\big[(\bI_k-\Omega(\kappa))^{-1}\big]_{ab}\tilde\bT_p(\kappa)\Lambda_b\tilde\bT_p(\kappa).
\]
Consequently,
\[
\frac1p\E\Tr\!\left(\tilde\bQ_p(\kappa)\Lambda_a\tilde\bQ_p(\kappa)\right)
=\sum_{b=1}^k\big[(\bI_k-\Omega(\kappa))^{-1}\big]_{ab}\frac1p\Tr\!\left(\tilde\bT_p(\kappa)\Lambda_b\tilde\bT_p(\kappa)\right)
+o(1).
\]
The same conclusions hold with $\tilde\bQ_p(\kappa)$
replaced by  $\tilde\bQ_p^{(-J)}(\kappa)$, where $J$ is a deterministic
set of indices of fixed cardinality.
\end{lemma}

The matrix deterministic equivalent follows from Proposition 1.6 of \cite{benaych2016spectral}, expressed in our notation.
The expectation statement follows from the uniform boundedness
of the normalized traces for fixed $\kappa>0$.

\begin{lemma}[Adapted from Lemma B.26 of \cite{bai2010spectral}]
\label{quad_form_lemma}
Let $X\in\mathbb R^p$ have independent, mean-zero,
unit-variance entries satisfying $\E|X_i|^k\le\nu_k$
for the relevant moment orders $k$. Then, for any
deterministic $A\in\mathbb R^{p\times p}$ and any
integer $k\ge1$ for which these moments exist,
\[
\E\bigl|X^\top AX-\tr(A)\bigr|^k
\le C_k\left[
\bigl(\nu_4\tr(AA^\top)\bigr)^{k/2}
+\nu_{2k}\bigl(\tr(AA^\top)\bigr)^{k/2}
\right],
\]
where $C_k$ is independent of $p$. In particular, if
$\|A\|\le1$ and the entries of $X$ have uniformly
bounded eighth moments, then
\[
\E\bigl|X^\top AX-\tr(A)\bigr|^4
\le Cp^2
\]
for a constant $C$ independent of $p$. Consequently,
when $p\to\infty$ through positive integers,
\[
\frac{X^\top AX}{p}-\frac{\tr(A)}{p}
\xrightarrow{\mathrm{a.s.}}0.
\]
\end{lemma}

\section{Additional proofs}

\subsection{Proof of Theorem~\ref{thm:joint_clt}}\label{apx:proof_join_clt}
\begin{proof}
Using Woodbury's matrix identity, one can easily check that for the ridge estimator we have
\begin{align}
    \widehat{\bbeta}
    =\frac1{\sqrt p}\bQ_p(\kappa)\PPhi^\top\Y_N
    =\frac1{\sqrt p}\PPhi^\top\tilde\bQ_p(\kappa)\Y_N.
\end{align}
Let $\{\pphi_l\}_{l=1}^p$ be the columns of $p^{-1/2}\PPhi$. In order to prove the theorem, we separate
the columns indexed by $J=\{j_1,\ldots,j_q\}$ from the rest. For this, we define
\begin{align}
    \tilde\bQ_p^{(-J)}(\kappa)
    &:= \left( \sum_{l\notin J}\pphi_l\pphi_l^\top+\kappa\bI_N    \right)^{-1},    &
    \Y_N^{(-J)}
    &:= \sum_{l\notin J}\pphi_l\beta_l+\Z_N.
\end{align}
We collect the selected columns and coordinates into
\begin{align}
     \PPhi^{(J)}&:=(\pphi_{j_1},\ldots,\pphi_{j_q}),
    &
    \bbeta_J &:=(\beta_{j_1},\ldots,\beta_{j_q})^\top,
    &
    \widehat{\bbeta}_J&:=(\hat\beta_{j_1},\ldots,\hat\beta_{j_q})^\top.
\end{align}

To express $\widehat{\bbeta}_J
= \PPhi^{(J)\top}\tilde\bQ_p(\kappa)\Y_N$
through these variables, we define $q\times q$ matrix
\begin{align}
    \mathbf H_J:=\left(\bI_q+ \PPhi^{(J)\top}\tilde\bQ_p^{(-J)}(\kappa) \PPhi^{(J)} \right)^{-1}.
\end{align}
By Woodbury's identity, we obtain
\begin{align}
    \tilde\bQ_p(\kappa)
    =\tilde\bQ_p^{(-J)}(\kappa)-\tilde\bQ_p^{(-J)}(\kappa) \PPhi^{(J)}\mathbf H_J \PPhi^{(J)\top}\tilde\bQ_p^{(-J)}(\kappa).
\end{align}
Multiplying on the left by $ \PPhi^{(J)\top}$ and using
\[
     \PPhi^{(J)\top}\tilde\bQ_p^{(-J)}(\kappa) \PPhi^{(J)}
    =\mathbf H_J^{-1}-\bI_q,
\]
we get
\begin{align}
     \PPhi^{(J)\top}\tilde\bQ_p(\kappa)
    =
    \mathbf H_J \PPhi^{(J)\top}\tilde\bQ_p^{(-J)}(\kappa).
\end{align}
Plugging this into $\widehat{\bbeta}_J$, we obtain
\begin{align}
    \widehat{\bbeta}_J
    =\mathbf H_J \PPhi^{(J)\top}\tilde\bQ_p^{(-J)}(\kappa)\left(\Y_N^{(-J)}+ \PPhi^{(J)}\bbeta_J\right)
=(\bI_q-\mathbf H_J)\bbeta_J+\mathbf H_J \PPhi^{(J)\top}\tilde\bQ_p^{(-J)}(\kappa)\Y_N^{(-J)}.
\end{align}
Define the diagonal matrices
\begin{align}
    \mathbf C_J
    &:=\operatorname{diag}
    \left(\mathcal C_{a_1},\ldots,\mathcal C_{a_q}\right),
    &
    \mathbf D_J
    &:=\operatorname{diag}
    \left(\kappa m_{a_1}(\kappa),\ldots,\kappa m_{a_q}(\kappa)\right).
\end{align}
Since $\mathbf C_J(\bI_q-\mathbf D_J)=\bI_q$, we obtain
\begin{align}\label{eq:prel_dec_J}
    \mathbf C_J\widehat{\bbeta}_J-\bbeta_J
    &=\mathbf C_J(\mathbf D_J-\mathbf H_J)\bbeta_J+\mathbf C_J\mathbf H_J\PPhi^{(J)\top}\tilde\bQ_p^{(-J)}(\kappa)\Y_N^{(-J)}.
\end{align}

In the next step we want to show that the first term on the r.h.s. is asymptotically small and substitute $\mathbf{H}_J$ with $\mathbf{D}_J$ in the second term. For this we use Woodbury's identity again, and obtain
\begin{align}
    \kappa\bQ_p(\kappa)
    =\bI_p-\frac1p\PPhi^\top\tilde\bQ_p(\kappa)\PPhi.
\end{align}
Taking the submatrix indexed by $J$ gives
\begin{align*}
    \kappa[\bQ_p(\kappa)]_{J,J}
    &=\bI_q-\PPhi^{(J)\top}\tilde\bQ_p(\kappa)\PPhi^{(J)}
    =\bI_q-\mathbf H_J\PPhi^{(J)\top}\tilde\bQ_p^{(-J)}(\kappa)\PPhi^{(J)}\\
    &    =\bI_q-\mathbf H_J(\mathbf H_J^{-1}-\bI_q)
    =\mathbf H_J.
\end{align*}
By the deterministic equivalent for $\bQ_p(\kappa)$, for every $r,s\in\{1,\ldots,q\}$,
\begin{align}
    [\bQ_p(\kappa)]_{j_r,j_s}
    =\delta_{rs}m_{a_r}(\kappa)+o_p(1).
\end{align}
Since $q$ is fixed, this implies
\begin{align}\label{eq:det_H_J}
    \|\mathbf H_J-\mathbf D_J\|=o_p(1).
\end{align}

To substitute this approximation into \eqref{eq:prel_dec_J},
we need to show that
\begin{align}
    \left\|\PPhi^{(J)\top}\tilde\bQ_p^{(-J)}(\kappa)\Y_N^{(-J)}\right\|=O_p(1).
\end{align}
Indeed, conditionally on
\[
    \mathcal F_{-J}
    :=\sigma\bigl(\{\pphi_l:l\notin J\},\Z_N\bigr),
\]
each coordinate
$\pphi_{j_r}^\top\tilde\bQ_p^{(-J)}(\kappa)\Y_N^{(-J)}$
is centered Gaussian, with conditional variance
\begin{align}
    \frac1p(\Y_N^{(-J)})^\top
    \tilde\bQ_p^{(-J)}(\kappa)\Lambda_{a_r}
    \tilde\bQ_p^{(-J)}(\kappa)\Y_N^{(-J)}
    \leq
    \frac{\|\Lambda_{a_r}\|}{p\kappa^2}
    \|\Y_N^{(-J)}\|^2.
\end{align}
The bounded covariance norms, $N/p=O(1)$ and
$\|\bbeta\|^2=O(p)$ imply
\[
    \mathbb E\|\Y_N^{(-J)}\|^2=O(p).
\]
Hence each coordinate has a bounded second moment.
Since $q$ is fixed, the required $O_p(1)$ bound follows.

Moreover, $\|\mathbf C_J\|=O(1)$ and
$\|\bbeta_J\|=O(1)$, by the boundedness assumptions on
the selected signal coordinates.
Thus our decomposition becomes
\begin{align}\label{eq:cbeta-beta_J}
    \mathbf C_J\widehat{\bbeta}_J-\bbeta_J
    =\mathbf C_J\mathbf D_J\PPhi^{(J)\top}\tilde\bQ_p^{(-J)}(\kappa)\Y_N^{(-J)}
    +o_p(1),
\end{align}
where the remainder converges to zero in probability
in Euclidean norm.

Conditionally on $\mathcal F_{-J}$, the vector $\mathbf C_J\mathbf D_J\PPhi^{(J)\top}\tilde\bQ_p^{(-J)}(\kappa)\Y_N^{(-J)}$
is centered Gaussian with diagonal conditional covariance matrix
\begin{align}
    \mathbf V_{J,p}(\kappa)
    :=\operatorname{diag}\left(V_{a_1,p}^{(-J)}(\kappa),\ldots,V_{a_q,p}^{(-J)}(\kappa)\right),
\end{align}
where
\begin{align}\label{eq:conditional_variance_J}
    V_{a,p}^{(-J)}(\kappa)
    :=\frac{\mathcal C_a^2\kappa^2m_a(\kappa)^2}{p}(\Y_N^{(-J)})^\top\tilde\bQ_p^{(-J)}(\kappa)\Lambda_a
    \tilde\bQ_p^{(-J)}(\kappa)\Y_N^{(-J)}.
\end{align}

It remains to prove that, for every $r=1,\ldots,q$,
\begin{align}\label{eq:variance_convergence_J}
    V_{a_r,p}^{(-J)}(\kappa)\xrightarrow{\mathbb P}V_{a_r}(\kappa).
\end{align}
Indeed, once \eqref{eq:variance_convergence_J} is established, we use L\'evy's  continuity theorem and check that for every $\mathbf t=(t_1,\ldots,t_q)^\top\in\mathbb R^q$,
\begin{align}
    \mathbb E\exp\left(i\mathbf t^\top\mathbf C_J\mathbf D_J\PPhi^{(J)\top}\tilde\bQ_p^{(-J)}(\kappa)\Y_N^{(-J)}\right)
    =\mathbb E\exp\left(-\frac12\sum_{r=1}^qt_r^2 V_{a_r,p}^{(-J)}(\kappa)\right)
\rightarrow
    \exp\left(-\frac12\sum_{r=1}^q t_r^2V_{a_r}(\kappa)\right).
\end{align}
Here we used convergence in probability of the conditional variances and boundedness of the exponential.
Thus
\[
  \mathbf C_J\mathbf D_J\PPhi^{(J)\top}\tilde\bQ_p^{(-J)}(\kappa)\Y_N^{(-J)}
    \xrightarrow{\mathrm d}\mathcal N\left(0,\operatorname{diag}\left(V_{a_1}(\kappa),\ldots,V_{a_q}(\kappa)\right)
    \right),
\]
and Slutsky's theorem handles the $o_p(1)$ remainder in \eqref{eq:cbeta-beta_J}.

In order to prove \eqref{eq:variance_convergence_J}, we first show that $V_{a,p}^{(-J)}(\kappa)$ concentrates, and then
that its expectation converges to $V_a(\kappa)$.

For this we use the following concentration result.
\begin{lemma}\label{lem:concent_J}
For fixed $\kappa>0$ and a deterministic set $J$ of fixed cardinality, under the assumptions of Theorem~\ref{thm:joint_clt},
\begin{align}
    \Var\left(V_{a,p}^{(-J)}(\kappa)\right)=O(p^{-1}),    \qquad a=1,\ldots,k.
\end{align}
\end{lemma}
For convenience, we postpone the proof to the Appendix~\ref{apx:proof_conc_lem}.

Let us now focus on $\mathbb E[V_{a,p}^{(-J)}(\kappa)]$. Plugging in the expression for $\Y_N^{(-J)}= \sum_{l\notin J}\pphi_l\beta_l+\Z_N$, we obtain
\begin{align}
    \mathbb E[V_{a,p}^{(-J)}(\kappa)]
    =    \mathcal C_a^2\kappa^2m_a(\kappa)^2    \left(A_{a,p}^{(-J)}+B_{a,p}^{(-J)}\right),
\end{align}
where
\begin{align}
    A_{a,p}^{(-J)}
    &:= \frac1p\sum_{l,m\notin J}\beta_l\beta_m
    \mathbb E\left[
        \pphi_l^\top\tilde\bQ_p^{(-J)}(\kappa)
        \Lambda_a\tilde\bQ_p^{(-J)}(\kappa)\pphi_m
    \right],\\
    B_{a,p}^{(-J)}
    &:=\frac{\sigma^2}{p}\mathbb E\Tr\left[
        \tilde\bQ_p^{(-J)}(\kappa)\Lambda_a
        \tilde\bQ_p^{(-J)}(\kappa)
    \right].
\end{align}
Here the mixed signal--noise terms vanish because  $\Z_N$ is centered and independent of the design,
and the noise quadratic form reduces to the trace because $\mathbb E[\Z_N\Z_N^\top]=\sigma^2\bI_N$.

We first consider $B_{a,p}^{(-J)}$. Removing a fixed number of columns does not change
the normalized double-resolvent trace by more than
$O(p^{-1})$. Indeed, $\tilde\bQ_p^{(-J)}(\kappa)-\tilde\bQ_p(\kappa)$
has rank at most $q$ and operator norm at most $\kappa^{-1}$. Consequently,
\begin{align}
    \frac1p\left|
    \Tr\left[\tilde\bQ_p^{(-J)}(\kappa)\Lambda_a\tilde\bQ_p^{(-J)}(\kappa)
        - \tilde\bQ_p(\kappa)\Lambda_a\tilde\bQ_p(\kappa) \right]\right|
    \leq\frac{2q\|\Lambda_a\|}{p\kappa^2}.
\end{align}
Applying Lemma~\ref{lem:double_resolvent}, we obtain
\begin{align}
    B_{a,p}^{(-J)}
    =\sigma^2\sum_{b=1}^k\big[(\bI_k-\Omega(\kappa))^{-1}\big]_{ab}
    \frac1p\Tr\left[\tilde\bT_p(\kappa)\Lambda_b\tilde\bT_p(\kappa)\right]+o(1).
\end{align}

Now let us look at $A_{a,p}^{(-J)}$. For $l\ne m$, changing $\pphi_l$ to $-\pphi_l$ leaves its distribution and $\tilde\bQ_p^{(-J)}(\kappa)$ unchanged, but changes the sign of the corresponding
bilinear form. Therefore, its expectation is zero, and all non-diagonal terms will vanish, i.e.
\begin{align}
    A_{a,p}^{(-J)}
    =\frac1p\sum_{l\notin J}\beta_l^2\mathbb E\left[\pphi_l^\top\tilde\bQ_p^{(-J)}(\kappa)
        \Lambda_a\tilde\bQ_p^{(-J)}(\kappa)\pphi_l\right].
\end{align}
For each block $b$, choose $i_b\in I_b\setminus J$.
By exchangeability of the remaining columns within each block,
\begin{align}
    A_{a,p}^{(-J)}
    =
    \frac1p\sum_{b=1}^k
    \left(
        \|\bbeta^{(b)}\|^2
        -\sum_{j\in J\cap I_b}\beta_j^2
    \right)
    \mathbb E\left[
        \pphi_{i_b}^\top\tilde\bQ_p^{(-J)}(\kappa)
        \Lambda_a\tilde\bQ_p^{(-J)}(\kappa)\pphi_{i_b}
    \right].
\end{align}
Since $q$ is fixed and the selected coefficients are
bounded, $\sum_{j\in J\cap I_b}\beta_j^2=O(1)$.
Moreover,
\[
    0\leq
    \mathbb E\left[
        \pphi_{i_b}^\top\tilde\bQ_p^{(-J)}(\kappa)
        \Lambda_a\tilde\bQ_p^{(-J)}(\kappa)\pphi_{i_b}
    \right]
    \leq
    \frac{\|\Lambda_a\|}{\kappa^2}
    \frac{\Tr\Lambda_b}{p}
    =O(1).
\]
Thus
\begin{align}\label{eq:signal_variance_J}
    A_{a,p}^{(-J)}
    =
    \frac1p\sum_{b=1}^k\|\bbeta^{(b)}\|^2
    \mathbb E\left[
        \pphi_{i_b}^\top\tilde\bQ_p^{(-J)}(\kappa)
        \Lambda_a\tilde\bQ_p^{(-J)}(\kappa)\pphi_{i_b}
    \right]+O(p^{-1}).
\end{align}
Since $\pphi_{i_b}$ and $\tilde\bQ_p^{(-J)}(\kappa)$ are not independent, in order to find expectation, we remove the column $i_b$ from the resolvent, and define
\begin{align}
    \tilde\bQ_p^{(-J\cup \{i_b\})}(\kappa)
    &:=
    \left(
        \sum_{l\notin J\cup\{i_b\}}\pphi_l\pphi_l^\top
        +\kappa\bI_N
    \right)^{-1},\\
    \eta_{J,i_b}
    &:=
    \left(
        1+\pphi_{i_b}^\top
        \tilde\bQ_p^{(-J\cup \{i_b\})}(\kappa)\pphi_{i_b}
    \right)^{-1}.
\end{align}
By the Sherman--Morrison identity,
\begin{align}
    \pphi_{i_b}^\top\tilde\bQ_p^{(-J)}(\kappa)
    =\eta_{J,i_b}\pphi_{i_b}^\top\tilde\bQ_p^{(-J\cup \{i_b\})}(\kappa).
\end{align}
Consequently, \eqref{eq:signal_variance_J} becomes
\begin{align}
    A_{a,p}^{(-J)}
    =
    \frac1p\sum_{b=1}^k\|\bbeta^{(b)}\|^2
    \mathbb E\left[
        \eta_{J,i_b}^2
        \pphi_{i_b}^\top\tilde\bQ_p^{(-J\cup \{i_b\})}(\kappa)
        \Lambda_a\tilde\bQ_p^{(-J\cup \{i_b\})}(\kappa)
        \pphi_{i_b}
    \right]+O(p^{-1}).
\end{align}
Arguing as in \eqref{eq:det_H_J}, now for the design
with the columns indexed by $J$ removed, we obtain
\[
    \eta_{J,i_b}=\kappa m_b(\kappa)+o_p(1).
\]
Here removing a fixed number of columns does not change
the deterministic equivalent at the $o_p(1)$ level.
To justify this replacement inside the expectation,
note that $0<\eta_{J,i_b}\leq1$ and
$0<\kappa m_b(\kappa)\leq1$. Thus
\[
    \mathbb E\left|
        \eta_{J,i_b}^2-\kappa^2m_b(\kappa)^2
    \right|^2=o(1).
\]
Moreover,
\[
    0\leq
    \pphi_{i_b}^\top\tilde\bQ_p^{(-J\cup \{i_b\})}(\kappa)
    \Lambda_a\tilde\bQ_p^{(-J\cup \{i_b\})}(\kappa)\pphi_{i_b}
    \leq
    \frac{\|\Lambda_a\|}{\kappa^2}\|\pphi_{i_b}\|^2,
\]
and $\mathbb E\|\pphi_{i_b}\|^4=O(1)$.
By Cauchy--Schwarz and $\|\bbeta\|^2/p=O(1)$,
we obtain
\begin{align}
    A_{a,p}^{(-J)}
    =
    \frac{\kappa^2}{p}
    \sum_{b=1}^k m_b(\kappa)^2\|\bbeta^{(b)}\|^2
    \mathbb E\left[
        \pphi_{i_b}^\top\tilde\bQ_p^{(-J\cup \{i_b\})}(\kappa)
        \Lambda_a\tilde\bQ_p^{(-J\cup \{i_b\})}(\kappa)\pphi_{i_b}
    \right]+o(1).
\end{align}
Finally, $\pphi_{i_b}$ is independent of
$\tilde\bQ_p^{(-J\cup \{i_b\})}(\kappa)$.
Taking the expectation with respect to $\pphi_{i_b}$ we get
\begin{align}\label{eq:signal_trace_J}
    A_{a,p}^{(-J)}
    =
    \frac{\kappa^2}{p^2}
    \sum_{b=1}^k m_b(\kappa)^2\|\bbeta^{(b)}\|^2
    \mathbb E\Tr\left[
        \Lambda_b\tilde\bQ_p^{(-J\cup \{i_b\})}(\kappa)
        \Lambda_a\tilde\bQ_p^{(-J\cup \{i_b\})}(\kappa)
    \right]+o(1).
\end{align}
To evaluate the remaining expectation, we apply
Lemma~\ref{lem:double_resolvent} with the deterministic
test matrix $\Lambda_b$. Removing $q+1$ columns changes
the normalized double-resolvent trace by $O(p^{-1})$.
Thus,
\begin{align}
    &\frac1p\mathbb E\Tr\left[
        \Lambda_b\tilde\bQ_p^{(-J\cup \{i_b\})}(\kappa)
        \Lambda_a\tilde\bQ_p^{(-J\cup \{i_b\})}(\kappa)
    \right]
=
    \sum_{d=1}^k
    \big[(\bI_k-\Omega(\kappa))^{-1}\big]_{ad}
    \frac1p\Tr\left[
        \Lambda_b\tilde\bT_p(\kappa)
        \Lambda_d\tilde\bT_p(\kappa)
    \right]+o(1).
\end{align}
Passing to expectations is justified by the uniform bound
\[
    0\leq
    \frac1p\Tr\left[
        \Lambda_b\tilde\bQ_p^{(-J\cup \{i_b\})}(\kappa)
        \Lambda_a\tilde\bQ_p^{(-J\cup \{i_b\})}(\kappa)
    \right]
    \leq
    \frac{N}{p\kappa^2}\|\Lambda_a\|\|\Lambda_b\|.
\]
Using the definition of $\Omega(\kappa)$, we can simplify
\[
    \frac1p\Tr\left[
        \Lambda_b\tilde\bT_p(\kappa)
        \Lambda_d\tilde\bT_p(\kappa)
    \right]
    =
    \frac{\Omega_{db}(\kappa)}
    {\ell_b\kappa^2m_b(\kappa)^2}.
\]
Consequently,
\begin{align}
    &\frac1p\mathbb E\Tr\left[
        \Lambda_b\tilde\bQ_p^{(-J\cup \{i_b\})}(\kappa)
        \Lambda_a\tilde\bQ_p^{(-J\cup \{i_b\})}(\kappa)
    \right]
=
    \frac{
        \big[(\bI_k-\Omega(\kappa))^{-1}
        \Omega(\kappa)\big]_{ab}
    }{\ell_b\kappa^2m_b(\kappa)^2}
    +o(1).
\end{align}
Substituting this into \eqref{eq:signal_trace_J}, we obtain
\begin{align}
    A_{a,p}^{(-J)}
    &=
    \sum_{b=1}^k
    \frac{\|\bbeta^{(b)}\|^2}{p\ell_b}
    \big[(\bI_k-\Omega(\kappa))^{-1}
    \Omega(\kappa)\big]_{ab}+o(1)
=
    \sum_{b=1}^k
    v_b\big[(\bI_k-\Omega(\kappa))^{-1}
    \Omega(\kappa)\big]_{ab}+o(1).
\end{align}
Finally, combining this with the expression for $B_{a,p}^{(-J)}$ gives
\begin{align}
    \mathbb E[V_{a,p}^{(-J)}(\kappa)]
    =
    \mathcal C_a^2\kappa^2m_a(\kappa)^2
    \Bigg[
    &\sum_{b=1}^k
    v_b\big[(\bI_k-\Omega(\kappa))^{-1}
    \Omega(\kappa)\big]_{ab}
    \nonumber\\
    &+
    \sigma^2\sum_{b=1}^k
    \big[(\bI_k-\Omega(\kappa))^{-1}\big]_{ab}
    \frac1p\Tr\left[
        \tilde\bT_p(\kappa)\Lambda_b\tilde\bT_p(\kappa)
    \right]
    \Bigg]+o(1).
\end{align}
Together with Lemma~\ref{lem:concent_J} and Chebyshev's inequality, we obtain
\[
    V_{a,p}^{(-J)}(\kappa)
    \xrightarrow{\mathbb P}V_a(\kappa).
\]
This proves \eqref{eq:variance_convergence_J}, and the
characteristic-function argument above completes the proof.
\end{proof}
\subsection{Proof of Lemma~\ref{det_eq_quad}}
\begin{proof}
 By the first-order deterministic equivalent
\eqref{eq:gram_det_eq}, for every fixed $t>0$,
\[
    \alpha^\intercal\bQ_p(t)\gamma
    =\sum_{a=1}^k
    \alpha^{(a)\intercal}\gamma^{(a)}m_a(t)+o_p(1).
\]
Applying this relation at $t=\kappa$ and $t=\kappa+h$, with $h>0$
fixed, subtracting and dividing by $h$, we obtain
\[
    \alpha^\intercal
    \frac{\bQ_p(\kappa)-\bQ_p(\kappa+h)}{h}\gamma
    =\sum_{a=1}^k \alpha^{(a)\intercal}\gamma^{(a)}
    \frac{m_a(\kappa)-m_a(\kappa+h)}{h}+o_p(1).
\]

To control the error as $h\downarrow0$, observe that for every $x\geq0$,
\[
    0\leq
    \frac{1}{(x+\kappa)^2}
    -\frac{1}{h}\left(\frac{1}{x+\kappa}
                     -\frac{1}{x+\kappa+h}\right)
    =\frac{h}{(x+\kappa)^2(x+\kappa+h)}
    \leq\frac{h}{\kappa^3}.
\]
Applying this bound to the eigenvalues of $\Phi_N^\intercal\Phi_N/p$
gives
\[
    \left\|\bQ_p(\kappa)^2
    -\frac{\bQ_p(\kappa)-\bQ_p(\kappa+h)}{h}\right\|
    \leq\frac{h}{\kappa^3}.
\]
Moreover, by \cite[Proposition~1.3]{benaych2016spectral}, in our
normalization $m_a(t)=\int_{[0,\infty)}(x+t)^{-1}\,\nu_{a,p}(dx)$
for a probability measure $\nu_{a,p}$. Integrating the same scalar
bound gives
\[
    \left|-m_a'(\kappa)
    -\frac{m_a(\kappa)-m_a(\kappa+h)}{h}\right|
    \leq\frac{h}{\kappa^3}.
\]
Since $\alpha$ and $\gamma$ are unit vectors, Cauchy--Schwarz gives
$\sum_{a=1}^k|\alpha^{(a)\intercal}\gamma^{(a)}|\leq1$.
Combining the preceding bounds, we obtain, for every fixed $h>0$,
\[
    \left|\alpha^\intercal\bQ_p(\kappa)^2\gamma
    +\sum_{a=1}^k\alpha^{(a)\intercal}\gamma^{(a)}m_a'(\kappa)\right|
    \leq\frac{2h}{\kappa^3}+o_p(1).
\]
Letting first $p\to\infty$ and then $h\downarrow0$ proves the first
assertion.

For the trace, we start from the first-order relation
\[
    \frac1p\Tr\bQ_p(t)
    =\sum_{a=1}^k\frac{|I_a|}{p}m_a(t)+o_p(1).
\]
The same argument applies, since the normalized trace of the resolvent
approximation error is bounded in absolute value by its operator norm
and the weights $|I_a|/p$ sum to one.
\end{proof}

\subsection{Proof of Lemma~\ref{lem:concent_J}}\label{apx:proof_conc_lem}
\begin{proof}
Here we prove the concentration of $V_{a,p}^{(-J)}(\kappa)$ for all $a=1\ldots,k$.
For convenience we recall that
\begin{align*}
    V_{a,p}^{(-J)}(\kappa)
    =    \frac{\mathcal C_a^2\kappa^2m_a(\kappa)^2}{p}(\Y_N^{(-J)})^\top\tilde\bQ_p^{(-J)}(\kappa)\Lambda_a\tilde\bQ_p^{(-J)}(\kappa)\Y_N^{(-J)},    \qquad
    \Y_N^{(-J)}=\sum_{l\notin J}\pphi_l\beta_l+\Z_N.
\end{align*} 
We prove the concentration in two steps: first with respect to the
noise $\Z_N$, and then with respect to the design columns
$\{\pphi_l:l\notin J\}$.
Denote by $\mathbb E_{\Z}$ expectation with respect to $\Z_N$,
keeping the design fixed. By the conditional variance decomposition,
\begin{align}\label{eq:variance_decomposition_J}
    \Var\bigl(V_{a,p}^{(-J)}(\kappa)\bigr)
    =\mathbb E\Big[\big(V_{a,p}^{(-J)}(\kappa)-\mathbb E_{\Z}V_{a,p}^{(-J)}(\kappa)\big)^2\Big]
+\mathbb E\Big[\big(\mathbb E_{\Z}V_{a,p}^{(-J)}(\kappa)-\mathbb EV_{a,p}^{(-J)}(\kappa)\big)^2\Big].
\end{align}
We will show that each term on the right-hand side is of order $O(p^{-1})$.

\textit{Concentration with respect to $\Z_N$.}
Since $\Z_N\sim\mathcal N(0,\sigma^2\bI_N)$ is independent of
the design, the Poincar\'e--Nash inequality gives
\begin{align}
    \mathbb E\Big[
        \big(V_{a,p}^{(-J)}(\kappa)
        -\mathbb E_{\Z}V_{a,p}^{(-J)}(\kappa)\big)^2
    \Big]
    \leq
    \sigma^2\sum_{\mu=1}^N
    \mathbb E\left[
        \left(\frac{\partial V_{a,p}^{(-J)}(\kappa)}
        {\partial Z_\mu}\right)^2
    \right].
\end{align}
The resolvent $\tilde\bQ_p^{(-J)}(\kappa)$ does not depend on
$\Z_N$, and $\partial\Y_N^{(-J)}/\partial Z_\mu=\mathbf e_\mu$.
Since $\tilde\bQ_p^{(-J)}(\kappa)\Lambda_a
\tilde\bQ_p^{(-J)}(\kappa)$ is symmetric, we obtain
\begin{align}
    \frac{\partial V_{a,p}^{(-J)}(\kappa)}{\partial Z_\mu}
    =
    \frac{2\mathcal C_a^2\kappa^2m_a(\kappa)^2}{p}
    \left(
        \tilde\bQ_p^{(-J)}(\kappa)\Lambda_a
        \tilde\bQ_p^{(-J)}(\kappa)\Y_N^{(-J)}
    \right)_\mu.
\end{align}
Summing the squared derivatives over $\mu$ and using
$\|\tilde\bQ_p^{(-J)}(\kappa)\|\leq\kappa^{-1}$ yields
\begin{align}
    \mathbb E\Big[\big(V_{a,p}^{(-J)}(\kappa)-\mathbb E_{\Z}V_{a,p}^{(-J)}(\kappa)\big)^2\Big]
&\leq\frac{4\sigma^2\mathcal C_a^4\kappa^4m_a(\kappa)^4}{p^2}\mathbb E\left\|\tilde\bQ_p^{(-J)}(\kappa)\Lambda_a\tilde\bQ_p^{(-J)}(\kappa)\Y_N^{(-J)}\right\|^2
\nonumber\\
    &\leq\frac{4\sigma^2\mathcal C_a^4m_a(\kappa)^4\|\Lambda_a\|^2}{p^2}\mathbb E\|\Y_N^{(-J)}\|^2.
\end{align}
Since the design columns are centered and mutually independent,
and $\Z_N$ is centered and independent of the design, we have
\begin{align*}
    \mathbb E\|\Y_N^{(-J)}\|^2
    &=\mathbb E\left[\left(\sum_{l\notin J}\pphi_l\beta_l+\Z_N\right)^\top\left(\sum_{m\notin J}\pphi_m\beta_m+\Z_N\right)
    \right]\\
    &=\sum_{l,m\notin J}\beta_l\beta_m\mathbb E[\pphi_l^\top\pphi_m]+2\sum_{l\notin J}\beta_l\mathbb E[\pphi_l^\top\Z_N]+\mathbb E[\Z_N^\top\Z_N]
    =\sum_{l\notin J}\beta_l^2\Tr\big(\mathbb E[\pphi_l\pphi_l^\top]\big)+\Tr\big(\mathbb E[\Z_N\Z_N^\top]\big).
\end{align*}
 For $l\in I_b$, we have
$\mathbb E[\pphi_l\pphi_l^\top]=p^{-1}\Lambda_b$,
while $\mathbb E[\Z_N\Z_N^\top]=\sigma^2\bI_N$. Therefore,
\begin{align*}
    \mathbb E\|\Y_N^{(-J)}\|^2
    =\frac1p\sum_{b=1}^k\sum_{l\in I_b\setminus J}\beta_l^2\Tr\Lambda_b+N\sigma^2
    &\leq \frac{N}{p}\sum_{b=1}^k\|\Lambda_b\|\sum_{l\in I_b\setminus J}\beta_l^2+N\sigma^2\\
    &\leq \frac{N}{p}\max_{1\leq b\leq k}\|\Lambda_b\|\|\bbeta\|^2+N\sigma^2
    =O(p),
\end{align*}
where we bounded $\Tr\Lambda_b\leq N\|\Lambda_b\|$ and
$\|\bbeta\|^2=O(p)$.
Thus, since the deterministic prefactor is bounded  for
fixed $\kappa>0$, we conclude that
\begin{align}\label{eq:noise_concentration_J}
    \mathbb E\Big[\big(V_{a,p}^{(-J)}(\kappa)-\mathbb E_{\Z}V_{a,p}^{(-J)}(\kappa)\big)^2\Big]
    =O(p^{-1}).
\end{align}

\textit{Concentration with respect to the design.}
For brevity, in this part we write
$\tilde\bQ=\tilde\bQ_p^{(-J)}(\kappa)$ and
$\Y=\Y_N^{(-J)}$, and denote
\[
    h:=\mathbb E_{\Z}V_{a,p}^{(-J)}(\kappa).
\]
The columns $\pphi_l$, $l\notin J$, are independent Gaussian
vectors with covariance $p^{-1}\Lambda_b$ for $l\in I_b$.
Thus, the Poincar\'e--Nash inequality gives
\begin{align}
    \Var(h)
    \leq\frac1p\sum_{b=1}^k\sum_{l\in I_b\setminus J}\mathbb E\left[(\nabla_{\pphi_l}h)^\top
        \Lambda_b\nabla_{\pphi_l}h\right]\leq\frac{\max_b\|\Lambda_b\|}{p}\sum_{l\notin J}\mathbb E\|\nabla_{\pphi_l}h\|^2.
\end{align}
Since, due to the Jensen inequality, $\left\|\nabla_{\pphi_l}\mathbb E_{\Z}(V_{a,p}^{(-J)}(\kappa))\right\|^2\leq\mathbb E_{\Z}\!\left[\left\|\nabla_{\pphi_l}V_{a,p}^{(-J)}(\kappa)\right\|^2
\right]$, we have
\begin{align}\label{eq:design_poincare_J}
    \Var(h)
    \leq\frac{\max_b\|\Lambda_b\|}{p}\sum_{l\notin J}\mathbb E\left\|\nabla_{\pphi_l}V_{a,p}^{(-J)}(\kappa)
    \right\|^2.
\end{align}
To compute these derivatives, recall that
\[
    \tilde\bQ
    =\left(\sum_{l\notin J}\pphi_l\pphi_l^\top
        +\kappa\bI_N\right)^{-1}.
\]
For $l\notin J$ and $\mu=1,\ldots,N$, the derivative of the
resolvent gives
\begin{align}
    \frac{\partial\tilde\bQ}{\partial\phi_{\mu l}}
    &=
    -\tilde\bQ\big(\mathbf e_\mu\pphi_l^\top+\pphi_l\mathbf e_\mu^\top\big)\tilde\bQ,
    &
    \frac{\partial\Y}{\partial\phi_{\mu l}}
    &=\beta_l\mathbf e_\mu,
\end{align}
where $\phi_{\mu l}$ denotes the $\mu$-th entry of $\pphi_l$. After straightforward calculations we obtain
\begin{align}
    \frac{\partial V_{a,p}^{(-J)}(\kappa)}{\partial\phi_{\mu l}}    
    =\frac{2\mathcal C_a^2\kappa^2m_a(\kappa)^2}{p}\Big[        \beta_l(\tilde\bQ\Lambda_a\tilde\bQ\Y)_\mu-(\pphi_l^\top\tilde\bQ\Y)(\tilde\bQ\Lambda_a\tilde\bQ\Y)_\mu
-(\pphi_l^\top\tilde\bQ\Lambda_a\tilde\bQ\Y)(\tilde\bQ\Y)_\mu\Big].
\end{align}
Consequently, using
$\|\mathbf u+\mathbf v+\mathbf w\|^2
\leq3(\|\mathbf u\|^2+\|\mathbf v\|^2+\|\mathbf w\|^2)$
in \eqref{eq:design_poincare_J}, we obtain
\begin{align}\label{eq:design_three_terms_J}
    \Var(h)
    \leq\frac{C}{p^3}\mathbb E\Big[\sum_{l\notin J}\beta_l^2\|\tilde\bQ\Lambda_a\tilde\bQ\Y\|^2
+\sum_{l\notin J}(\pphi_l^\top\tilde\bQ\Y)^2\|\tilde\bQ\Lambda_a\tilde\bQ\Y\|^2+\sum_{l\notin J}
            (\pphi_l^\top\tilde\bQ\Lambda_a\tilde\bQ\Y)^2\|\tilde\bQ\Y\|^2\Big],
\end{align}
where $C$ is independent of $p$.

For the first term, as many times before, we use $\|\tilde\bQ\|\leq\kappa^{-1}$,
$\|\bbeta\|^2=O(p)$, and $\mathbb E\|\Y\|^2=O(p)$, this gives us
\begin{align}
    \mathbb E\left[\sum_{l\notin J}\beta_l^2\|\tilde\bQ\Lambda_a\tilde\bQ\Y\|^2\right]
    \leq\frac{\|\Lambda_a\|^2}{\kappa^4}\|\bbeta\|^2\mathbb E\|\Y\|^2
    =O(p^2).
\end{align}

Using $\sum_{l\notin J}\pphi_l\pphi_l^\top=\tilde\bQ^{-1}-\kappa\bI_N$, we rewrite the second term as
\begin{align}
   \mathbb E &\left[\sum_{l\notin J}(\pphi_l^\top\tilde\bQ\Y)^2\|\tilde\bQ\Lambda_a\tilde\bQ\Y\|^2\right]
  =\mathbb E\left[\Y^\top(\tilde\bQ-\kappa\tilde\bQ^2)\Y\|\tilde\bQ\Lambda_a\tilde\bQ\Y\|^2\right]
    \nonumber\\
    &\quad\leq\mathbb E\left[\|\tilde\bQ-\kappa\tilde\bQ^2\|\|\tilde\bQ\|^4\|\Lambda_a\|^2\|\Y\|^4\right]
\leq\|\Lambda_a\|^2\mathbb E\left[\|\tilde\bQ\|^5\|\Y\|^4\right]
\leq\frac{\|\Lambda_a\|^2}{\kappa^5}\mathbb E\|\Y\|^4,
\end{align}
where we used
$0\preceq\tilde\bQ-\kappa\tilde\bQ^2\preceq\tilde\bQ$
and $\|\tilde\bQ\|\leq\kappa^{-1}$.
Similarly, for the third term,
\begin{align}
   \mathbb E\left[\sum_{l\notin J}(\pphi_l^\top\tilde\bQ\Lambda_a\tilde\bQ\Y)^2\|\tilde\bQ\Y\|^2\right]
 &=\mathbb E\left[\Y^\top\tilde\bQ\Lambda_a(\tilde\bQ-\kappa\tilde\bQ^2)\Lambda_a\tilde\bQ\Y\|\tilde\bQ\Y\|^2
    \right]\\\nonumber
&\leq\mathbb E\left[\|\tilde\bQ-\kappa\tilde\bQ^2\|\|\tilde\bQ\|^4\|\Lambda_a\|^2\|\Y\|^4\right]
\leq\|\Lambda_a\|^2\mathbb E\left[\|\tilde\bQ\|^5\|\Y\|^4\right]
    \leq\frac{\|\Lambda_a\|^2}{\kappa^5}\mathbb E\|\Y\|^4.
\end{align}
Since $\Y$ is a centered Gaussian vector, we conclude $\mathbb E\|\Y\|^4=O(p^2).$

Therefore all three terms inside the expectation in
\eqref{eq:design_three_terms_J} are of order $O(p^2)$, and
\begin{align}
    \Var(h)=O(p^{-1}).
\end{align}
Combining this estimate with \eqref{eq:noise_concentration_J} and \eqref{eq:variance_decomposition_J} proves
\[
    \Var\bigl(V_{a,p}^{(-J)}(\kappa)\bigr)=O(p^{-1}).
\]
\end{proof}

\section{Further results on Bayesian estimation}\label{app:bayes_additional}

We derive the Gaussian formulas for the measurement MMSE and mutual
information, relate the variance in the joint CLT to the Bayesian
fixed-point equations, and give an upper bound on the blockwise Bayes
risk for non-Gaussian priors.

\subsection{Measurement MMSE and mutual information}

The parameter MMSE is normalized by $p$, whereas the measurement MMSE
is normalized by $N$. With
\[
 \operatorname{ymmse}_p
 :=\frac1N\mathbb E\left\|
 X_p\{\bbeta-\mathbb E(\bbeta\mid\Y_N,\PPhi)\}\right\|^2,
\]
the exact posterior covariance yields
\begin{align}\label{eq:bayes_conn_ymmse_trace}
 \operatorname{ymmse}_p
 &=\frac{\sigma^2}{N}\mathbb E\Tr\{S_pQ_p(\kappa_{\mathrm B})\}
  =\frac{\sigma^2}{N}
       \left[p-\kappa_{\mathrm B}\mathbb E\Tr Q_p(\kappa_{\mathrm B})\right]
       \nonumber\\
 &\rightarrow
 \frac{\sigma^2}{c}
       \left(1-\frac1\rho\sum_a\ell_ar_{2,a}\right).
\end{align}
Writing $A(\theta):=\sum_a\ell_af_a(\theta)r_{2,a}$,
multiplication of the second equation in
\eqref{eq:bayes_conn_bayes_fp} by $\ell_ar_{2,a}$ and summation give
\[
 1-\frac1\rho\sum_a\ell_ar_{2,a}
 =\sum_a\ell_ar_{1,a}r_{2,a}
 =\frac c\pi\int_0^\pi
       \frac{A(\theta)}{\sigma^2+A(\theta)}\,d\theta.
\]
Consequently,
\begin{equation}\label{eq:bayes_conn_ymmse_integral}
 \lim_{p\to\infty}\operatorname{ymmse}_p
 =\frac{\sigma^2}{\pi}\int_0^\pi
       \frac{\sum_a\ell_af_a(\theta)r_{2,a}}
            {\sigma^2+\sum_a\ell_af_a(\theta)r_{2,a}}\,d\theta,
\end{equation}
which is the Gaussian specialization of
\cite[Eq.~(II.13)]{mutinf_time}. The derivation is direct and does
not infer $k$ block errors from a single measurement-error identity.

The correspondence also extends to the value of the normalized mutual
information. Here conditional mutual information includes averaging over
the random design. Gaussian conjugacy gives
\begin{equation}\label{eq:bayes_conn_gaussian_MI_exact}
 i_p:=\frac1p I(\bbeta;\Y_N\mid\PPhi)
 =\frac1{2p}\mathbb E\log\det
       \left(\bI_p+\frac{\rho}{\sigma^2}S_p\right).
\end{equation}

\begin{proposition}[Gaussian mutual-information limit]
\label{prop:bayes_conn_MI}
At the solution in \eqref{eq:bayes_conn_dictionary},
\begin{equation}\label{eq:bayes_conn_MI_limit}
 \begin{split}
 \lim_{p\to\infty}i_p
 ={}&\frac c{2\pi}\int_0^\pi
       \log\left(1+
          \frac{\sum_a\ell_af_a(\theta)r_{2,a}}{\sigma^2}
           \right)\,d\theta\\
 &-\frac12\sum_a\ell_ar_{1,a}r_{2,a}
  +\frac12\sum_a\ell_a\log(1+\rho r_{1,a}).
 \end{split}
\end{equation}
This is the Gaussian value of the potential in
\cite[Eq.~(II.4)]{mutinf_time} at its unique positive stationary point.
\end{proposition}

\begin{proof}
For $t>0$, put
\[
 J_p(t):=\frac1{2p}\mathbb E\log\det(\bI_p+tS_p),\qquad
 h_a(t):=\frac{\bar m_a(1/t)}t,
 \qquad
 D_t(\theta):=1+t\sum_a\ell_ah_a(t)f_a(\theta).
\]
The fixed-point equations are equivalent to
\[
 h_a(t)^{-1}=1+t\frac c\pi\int_0^\pi
                         \frac{f_a(\theta)}{D_t(\theta)}\,d\theta.
\]
The functions $h_a$ are differentiable for $t>0$, by stability of
the fixed-point system, and $h_a(t)\to1$ as $t\downarrow0$.
Differentiating the finite-dimensional log determinant and applying
the trace deterministic equivalent at $\kappa=1/t$ give
\[
 J_p'(t)=\frac1{2p}\mathbb E\Tr\{S_p(\bI_p+tS_p)^{-1}\}
 \rightarrow j(t):=\frac1{2t}\left(1-\sum_a\ell_ah_a(t)\right).
\]
Moreover,
\[
 0\le J_p'(t)\le\frac1{2p}\mathbb E\Tr S_p
  =\frac1{2p^2}\sum_a p_a\Tr\Lambda_a\le C,
\]
uniformly in $t\ge0$ and large $p$. Since $J_p(0)=0$, dominated
convergence yields $J_p(t)\to\int_0^t j(u)\,du$.

Define
\[
 F(t):=\frac c{2\pi}\int_0^\pi\log D_t(\theta)\,d\theta
       -\frac12\sum_a\ell_a(1-h_a(t))
       -\frac12\sum_a\ell_a\log h_a(t).
\]
The partial derivative of this expression with respect to $h_a$,
keeping $t$ fixed, equals
\[
 \frac{\ell_a}{2}\left[
   t\frac c\pi\int_0^\pi\frac{f_a(\theta)}{D_t(\theta)}\,d\theta
       +1-h_a(t)^{-1}\right]=0.
\]
Consequently the terms involving $h_a'(t)$ cancel, and
\[
 F'(t)=\frac c{2\pi}\int_0^\pi
       \frac{\sum_a\ell_ah_a(t)f_a(\theta)}{D_t(\theta)}\,d\theta
       =\frac1{2t}\sum_a\ell_a(1-h_a(t))=j(t).
\]
Also $F(t)\to0$ as $t\downarrow0$. Thus $J_p(t)\to F(t)$.
At $t=\rho/\sigma^2$ one has
$h_a=r_{2,a}/\rho$, $1-h_a=r_{1,a}r_{2,a}$, and
$-\log h_a=\log(1+\rho r_{1,a})$.
Substitution proves \eqref{eq:bayes_conn_MI_limit}.
This argument obtains the Gaussian information formula from the
resolvent equations without applying a compact-support hypothesis
to a Gaussian prior.
\end{proof}

\subsection{Bayesian errors versus the debiased CLT variance.}

Next we show that under the  further assumptions, the limiting variance of the debiased estimators can be related to the Bayesian fixed-point system.

\begin{proposition}[Matched-penalty simplification of the CLT variance]
\label{prop:bayes_conn_clt_variance}
If the block energies in Theorem~\ref{thm:joint_clt} satisfy
$v_a=\rho$ for every $a$, then, at $\kappa=\kappa_{\mathrm B}$,
\begin{equation}\label{eq:bayes_conn_V_identity}
 \quad
 V_a(\kappa_{\mathrm B})
 =\frac{\sigma^2\bar m_a(\kappa_{\mathrm B})}
        {1-\kappa_{\mathrm B}\bar m_a(\kappa_{\mathrm B})}
 =\frac{r_{2,a}}{1-r_{2,a}/\rho}
 =\frac1{r_{1,a}}.
 \quad
\end{equation}
In particular, if
$\overline{\mathcal C}_a:=(1-\kappa_{\mathrm B}\bar m_a)^{-1}$,
then $V_a=\overline{\mathcal C}_a r_{2,a}>r_{2,a}$.
The identity also holds for the CLT variance of a deterministic signal
with these block energies; Gaussianity of that signal is not needed
for this algebraic simplification.
\end{proposition}

\begin{proof}
All quantities in this proof are limiting quantities at
$\kappa=\kappa_{\mathrm B}$. Abbreviate
$m_a:=\bar m_a(\kappa)$, $\Omega:=\bar\Omega(\kappa)$,
$B:=(\bI_k-\Omega)^{-1}$, and set
\[
 d_a:=m_a^{-1}-\kappa,\qquad
 w_a:=m_a^{-1},\qquad
 \zeta_a:=\frac{c}{\pi\kappa^2}
        \int_0^\pi\frac{f_a(\theta)}{\tau(\theta,\kappa)^2}\,d\theta.
\]
Using the definition of $\Omega$ and
$\tau-1=\sum_b\ell_bm_bf_b$, we obtain
\[
 (\Omega\boldsymbol w)_a
 =\frac c\pi\int_0^\pi
     \frac{f_a(\theta)(\tau(\theta,\kappa)-1)}
          {\tau(\theta,\kappa)^2}\,d\theta
 =d_a-\kappa^2\zeta_a.
\]
Since $\boldsymbol d=\boldsymbol w-\kappa\boldsymbol1$,
\[
 (\bI_k-\Omega)\boldsymbol d
   =\kappa\Omega\boldsymbol1+\kappa^2\boldsymbol\zeta.
\]
At the matched penalty $\sigma^2=\rho\kappa$, so
\[
 B\{\rho\Omega\boldsymbol1+\sigma^2\boldsymbol\zeta\}
 =\frac\rho\kappa\boldsymbol d.
\]
The variance formula \eqref{eq:coordinate_variance} therefore reduces to
\[
 V_a=\overline{\mathcal C}_a^{\,2}\kappa^2m_a^2
                    \frac\rho\kappa d_a
     =\overline{\mathcal C}_a^{\,2}\sigma^2m_a(1-\kappa m_a)
     =\frac{\sigma^2m_a}{1-\kappa m_a}.
\]
Finally, \eqref{eq:bayes_conn_dictionary} gives
$d_a=\sigma^2r_{1,a}$ and $\sigma^2m_a=r_{2,a}$, proving
all expressions in \eqref{eq:bayes_conn_V_identity}.
\end{proof}

\begin{remark}

We note that the matrix $\Omega$ appearing in the CLT has a Bayesian interpretation. Eliminate $r_2$ from the Gaussian Bayesian system and define the
precision update map
\[
 \mathcal F_a(s)
 :=\frac c\pi\int_0^\pi
 \frac{f_a(\theta)}
 {\sigma^2+\sum_b\ell_bf_b(\theta)(\rho^{-1}+s_b)^{-1}}
 \,d\theta.
\]
At its fixed point $s=r_1$, direct differentiation gives
\begin{equation}\label{eq:bayes_conn_jacobian}
 \frac{\partial\mathcal F_a}{\partial s_b}(r_1)
 =\frac{c\ell_b r_{2,b}^2}{\pi}
  \int_0^\pi\frac{f_a(\theta)f_b(\theta)}
   {\{\sigma^2+\sum_d\ell_df_d(\theta)r_{2,d}\}^2}\,d\theta
 =\bar\Omega_{ab}(\kappa_{\mathrm B}).
\end{equation}
Thus the matrix appearing in the frequentist CLT is also the Jacobian
of the Gaussian Bayesian precision equations. The bound
$\rho(\bar\Omega)<1$ established in
Proposition~\ref{prop:explicit_variance_limit} implies local stability
of this fixed point. The inverse $(\bI_k-\bar\Omega)^{-1}$ describes
the linear response of the same equations to perturbations. This
identification concerns the fixed-point system; it is not a convergence
theorem for a message-passing algorithm.
\end{remark}

\subsection{Scope of the conclusion and non-Gaussian priors.}

The results above provide a direct proof of the Gaussian specialization
of the block-MMSE formula. The compact-support assumption in
\cite[Theorem~2.1]{mutinf_time} does not include a nondegenerate
Gaussian prior; Propositions~\ref{prop:bayes_conn_posterior}--\ref{prop:bayes_conn_MI}
therefore establish the Gaussian identities independently of that
assumption. They do not establish the block-MMSE conjecture for a
general non-Gaussian prior.

More precisely, suppose only that the prior is centered, independent of
the design, and has covariance $\rho\bI_p$. The conditional covariance
calculation still shows that ridge at $\kappa_{\mathrm B}$ is the
optimal estimator among those linear in $\Y_N$, with conditional
error covariance $\sigma^2Q_p(\kappa_{\mathrm B})$.
For example, this follows by minimizing the conditional quadratic risk
over the coefficient matrix of a linear estimator, using
$\operatorname{Cov}(\bbeta,\Y_N\mid\PPhi)=\rho X_p^\top$ and
$\operatorname{Cov}(\Y_N\mid\PPhi)=\rho X_pX_p^\top+\sigma^2\bI_N$.
Writing $\eta_p:=\mathbb E(\bbeta\mid\Y_N,\PPhi)$, conditional
expectation orthogonality gives the exact gap identity
\begin{equation}\label{eq:bayes_conn_nongaussian_gap}
 \frac{\sigma^2}{p_a}\mathbb E\Tr\{P_aQ_p(\kappa_{\mathrm B})\}
 =\frac1{p_a}\mathbb E\|P_a(\bbeta-\eta_p)\|^2
  +\frac1{p_a}\mathbb E
       \|P_a(\eta_p-\widehat{\bbeta}_{\kappa_{\mathrm B}})\|^2.
\end{equation}
Consequently the Gaussian-prior value gives an upper bound on the
non-Gaussian Bayes risk,
\[
 \limsup_{p\to\infty}\frac1{p_a}
       \mathbb E\|P_a(\bbeta-\eta_p)\|^2
 \le\sigma^2\bar m_a(\sigma^2/\rho).
\]
For a non-Gaussian prior $P_0^{\otimes p}$, the first equation in
\eqref{eq:bayes_conn_bayes_fp} retains its prior-dependent scalar
MMSE function, and the last term of
\eqref{eq:bayes_conn_nongaussian_gap} need not vanish.

\section{Additional results}

\begin{proposition}\label{prop:explicit_variance_limit}
For $\theta\in[-\pi,\pi]$, define
\[
    f_b(\theta)
    :=\frac{1}{1-2\lambda_b\cos\theta+\lambda_b^2}.
\]
For fixed $\kappa>0$, let
$\bar{m}_1(\kappa),\ldots,\bar m_k(\kappa)$ be the unique
positive solution of
\begin{align}\label{eq:limiting_fixed_point}
    \bar m_a(\kappa)
    &=
    \left(
        \kappa+
        \frac{c}{2\pi}\int_{-\pi}^{\pi}
        \frac{f_a(\theta)}{\tau(\theta,\kappa)}\,d\theta
    \right)^{-1},
    \qquad a=1,\ldots,k,
    \\
    \tau(\theta,\kappa)
    &:=
    1+\sum_{b=1}^k\ell_b\bar m_b(\kappa)f_b(\theta).
\end{align}
Then $m_a(\kappa)\to\bar m_a(\kappa)$ as $N,p\to\infty$.

Moreover, define the matrix $\bar\Omega(\kappa)$ with elements
\begin{align}\label{eq:limiting_Omega}
    \bar\Omega_{ab}(\kappa)
    :=
    \frac{c\ell_b\bar m_b(\kappa)^2}{2\pi}
    \int_{-\pi}^{\pi}
    \frac{f_a(\theta)f_b(\theta)}
         {\tau(\theta,\kappa)^2}\,d\theta.
\end{align}
The spectral radius of $\bar\Omega(\kappa)$ is strictly less than one, and the limit in \eqref{eq:coordinate_variance} exists and is given by
\begin{align}\label{eq:explicit_limiting_variance}
    V_a(\kappa)
    =\frac{\kappa^2\bar m_a(\kappa)^2}{(1-\kappa\bar m_a(\kappa))^2}\Bigg[&\sum_{b=1}^kv_b\big[
            (\bI_k-\bar\Omega(\kappa))^{-1}\bar\Omega(\kappa)\big]_{ab}+\frac{c\sigma^2}{2\pi\kappa^2}
        \sum_{b=1}^k\big[(\bI_k-\bar\Omega(\kappa))^{-1}\big]_{ab}\int_{-\pi}^{\pi}\frac{f_b(\theta)}
             {\tau(\theta,\kappa)^2}\,d\theta\Bigg].
\end{align}
\end{proposition}
\begin{proof}
We use the Toeplitz matrix results recalled in
\cite[Appendix~B]{mutinf_time}.
In particular, the covariance matrices $\Lambda_b$ are scaled
Kac--Murdock--Szeg\"o matrices with generating functions
\[
    f_b(\theta)
    =\frac{1}{1-2\lambda_b\cos\theta+\lambda_b^2}.
\]
The asymptotic trace formulas for products and inverses of
these matrices allow us to express the required limits as
integrals of their generating functions.
We first establish the convergence of $m_a(\kappa)$.

To make the dependence on the dimension explicit, write
$m_{a,p}(\kappa)$ for $m_a(\kappa)$ in this proof.
The finite-dimensional fixed-point equations give
\begin{align}
    m_{a,p}(\kappa)
    =
    \left(
        \kappa+\frac1p\Tr\left[
            \Lambda_a
            \left(
                \bI_N+\sum_{b=1}^k
                \frac{|I_b|}{p}m_{b,p}(\kappa)\Lambda_b
            \right)^{-1}
        \right]
    \right)^{-1}.
\end{align}
Since the matrix inside the inverse is bounded below by $\bI_N$, we have
\[
    \frac{1}{\kappa+(N/p)\|\Lambda_a\|}
    \leq m_{a,p}(\kappa)\leq\kappa^{-1}.
\]
Hence $(m_{a,p}(\kappa))_{a=1}^k$ remains in a compact subset of $(0,\infty)^k$. Along any convergent subsequence, the Toeplitz trace formulas recalled in \cite[Appendix~B]{mutinf_time}, together with the resolvent identity, allow us to pass to the limit in the fixed-point equations. Every subsequential limit therefore solves \eqref{eq:limiting_fixed_point}. By uniqueness of the positive solution of this canonical system, we conclude that
\[
    m_{a,p}(\kappa)\rightarrow\bar m_a(\kappa),
    \qquad a=1,\ldots,k.
\]

Using this convergence and the same Toeplitz trace formulas,
we obtain
\begin{align}
    \frac1p\Tr\left[
        \Lambda_a\tilde\bT_p(\kappa)
        \Lambda_b\tilde\bT_p(\kappa)
    \right]
    &\rightarrow
    \frac{c}{2\pi\kappa^2}
    \int_{-\pi}^{\pi}
    \frac{f_a(\theta)f_b(\theta)}
         {\tau(\theta,\kappa)^2}\,d\theta,
    \\
    \frac1p\Tr\left[
        \tilde\bT_p(\kappa)\Lambda_b\tilde\bT_p(\kappa)
    \right]
    &\rightarrow
    \frac{c}{2\pi\kappa^2}
    \int_{-\pi}^{\pi}
    \frac{f_b(\theta)}
         {\tau(\theta,\kappa)^2}\,d\theta.
\end{align}
In particular, $\Omega(\kappa)\to\bar\Omega(\kappa)$
entrywise.

To show that $\bI_k-\bar\Omega(\kappa)$ is invertible, set
$u_a=1/\bar m_a(\kappa)$ and
$\bU=\operatorname{diag}(u_1,\ldots,u_k)$. Then
\begin{align}
    [\bar\Omega(\kappa)\mathbf u]_a
    &=
    \frac{c}{2\pi}\int_{-\pi}^{\pi}
    \frac{f_a(\theta)(\tau(\theta,\kappa)-1)}
         {\tau(\theta,\kappa)^2}\,d\theta
    \nonumber\\
    &<
    \frac{c}{2\pi}\int_{-\pi}^{\pi}
    \frac{f_a(\theta)}{\tau(\theta,\kappa)}\,d\theta
    =u_a-\kappa<u_a,
\end{align}
where we used \eqref{eq:limiting_fixed_point}.
Since $\bar\Omega(\kappa)$ has nonnegative entries,
\[
    \rho(\bar\Omega(\kappa))
    =\rho(\bU^{-1}\bar\Omega(\kappa)\bU)
    \leq\|\bU^{-1}\bar\Omega(\kappa)\bU\|_\infty
    =\max_a\frac{[\bar\Omega(\kappa)\mathbf u]_a}{u_a}
    <1,
\]
where $\rho(\cdot)$ denotes the spectral radius and
$\|\cdot\|_\infty$ denotes the maximum absolute row sum.
Consequently, by continuity of matrix inversion,
\[
    (\bI_k-\Omega(\kappa))^{-1}
    \rightarrow
    (\bI_k-\bar\Omega(\kappa))^{-1}.
\]
Moreover, \eqref{eq:limiting_fixed_point} implies
$0<\kappa\bar m_a(\kappa)<1$, and hence
\[
    \mathcal C_a
    =\frac1{1-\kappa m_a(\kappa)}
    \rightarrow
    \frac1{1-\kappa\bar m_a(\kappa)}.
\]
Substituting these limits into \eqref{eq:coordinate_variance}
gives \eqref{eq:explicit_limiting_variance}.
\end{proof}

\end{document}